\documentclass[12pt]{article}

\usepackage{graphicx}    % For including images
\usepackage{amsmath} 
\usepackage{xcolor}
\usepackage{amssymb, amsthm} 
\usepackage{booktabs}% Additional math symbols
\usepackage{hyperref}    % For hyperlinks
\usepackage{geometry}    % Adjust margins
\usepackage{caption}     % For captions in figures and tables
\usepackage{cite}
\usepackage{algorithm}
\usepackage{algorithmicx}
\usepackage{algpseudocode}
\usepackage{pgfplots}
\usepackage{bm}
\usepackage{float}
\usetikzlibrary{patterns,positioning}
\usepackage{tikz}
\usetikzlibrary{decorations.markings}
\usetikzlibrary{matrix}
\usepackage{pgfplots}
\usetikzlibrary{shapes, arrows.meta, patterns}
\usepackage{mathrsfs}

\usepackage{subcaption}

\newcommand{\bfa}[1]{\boldsymbol{#1}} 			
			
\newcommand{\bfeps}{\boldsymbol{\epsilon}}

\newcommand{\bfsigma}{\boldsymbol{\sigma}}

\DeclareMathAlphabet{\mathpzc}{OT1}{pzc}{m}{it}

\DeclareMathOperator*{\argmin}{arg\, min}
\newcommand{\Th}{\mathscr{T}_{h}}

\newcommand{\Fh}{\mathscr{F}_{h}}

\newcommand{\refface}{\hat{\tau}} % Or \hat{K}, \hat{T}
\newcommand{\refcoord}{\hat{\mathbf{z}}}
\newcommand{\nodes}{\mathcal{I}_h}
\newcommand{\nodecoord}{\mathbf{z}}
\newcommand{\basis}{\zeta}

\newcommand{\nodesDboundary}{\mathcal{I}_{D,h}}
\newcommand{\FESpace}{\mathcal{S}_h}
\newcommand{\FESpaceHomD}{\mathcal{S}_{d,h}}
\newcommand{\FESpaceInhomD}{\mathcal{S}^n_{g,h}}
\newcommand{\FESpaceCrack}{\mathcal{S}^n_{cr,h}}
\newcommand{\FESpaceCrackSimple}{\mathcal{S}_{cr,h}}
\newcommand{\FESpaceInhomDSimple}{\mathcal{S}_{g,h}}
\newcommand{\FacesCR}{\mathscr{F}_h^{CR}}
\newcommand{\CrackRegion}{C_h}
\newcommand{\tolCR}{\Xi_{CR}}
\newcommand{\phasefield}{\psi_h}

\usepackage{authblk}

\usepackage{etoolbox}

\newcommand{\Hspace}[1]{\mathrm{H}^{#1}}
\newcommand{\Lspace}[1]{\mathrm{L}^{#1}}

\newtheorem*{dwf}{Discrete formulation}

\title{$\xi$-Based adaptive phase field model for quasi-static anti-plane fracture}
\author[1]{Maria P. Fernando}
\author[1,*]{S. M. Mallikarjunaiah}

\affil[1]{Department of Mathematics \& Statistics, Texas A\&M University-Corpus Christi, TX- 78412, USA}
\affil[*]{Corresponding author}
\affil[ ]{\textit{E-mail addresses:} \texttt{fpieo@islander.tamucc.edu} (M.P. Fernando), \texttt{M.Muddamallappa@tamucc.edu} (S.M. Mallikarjunaiah)}

\date{}

\begin{document}

\maketitle  

\begin{abstract}
The $\xi$-based spatially adaptive three-field variable phase-field model for quasi-static anti-plane crack propagation is introduced. A dynamically optimized regularization length is integrated to improve computational efficiency and accuracy in numerical approximations. A local adaptive mesh refinement strategy is developed, which maintains an optimal balance between mesh resolution and the accurate depiction of fractures using the \textsf{AT1} diffuse interface model. The total energy functional is comprised of three components: strain energy, surface energy, and a third term reliant on the damage zone's regularization length. The governing partial differential equations for mechanics and phase-field variables, derived from Euler-Lagrange, are discretized via the finite-element method. Two parameters functioning as penalty variables are incorporated; both are asymptotically estimated from the gradient of the phase-field variable. By these estimated parameters, mesh adaptivity is enhanced, ensuring the convergence of the numerical solution. Standard phase-field methods are shown by numerical results to be surpassed by the adaptive model; an accurate representation of fractures is provided, and computational costs are significantly lowered.  By employing the proposed spatially adaptive approach, a vastly larger regularization length parameter is achieved compared to other methods throughout the entire computation. 
\end{abstract}

\vspace{.1in}


\section{Introduction}

The study of cracks and fractures in materials, particularly within elastic and porous media, represents a significant and active area of research across applied mathematics and engineering \cite{gou2015modeling,ghosh2025finite,manohar2025convergence,rajagopal2011modeling,Mallikarjunaiah2015,ferguson2015,yoon2024finite,yoon2022MAM,yoon2022CNSNS,silling2000reformulation,vasilyeva2024generalized,manohar2024hp}. This field, known as \textit{fracture mechanics}, is crucial for ensuring the safety and reliability of everything from everyday products to massive infrastructure. By understanding how cracks form and grow, engineers can prevent catastrophic failures in structures like bridges, airplanes, and pipelines. This knowledge enables the design of more durable materials and components and helps determine the safe operating life of a structure, thereby preventing potential loss of life and economic devastation.

Modeling the behavior of cracks is an inherently challenging task due to the complexities involved in describing the entire process with a set of equations. The study of crack evolution has been approached using various techniques, including methods that explicitly define the crack-tip movement, such as the \textit{extended finite element method} \cite{fries2010extended} and the \textit{generalized finite element method}  \cite{babuvska2012stable}, among many others. Recently, \textit{variational regularization approaches}, more commonly known as \textit{phase-field approximation methods}, have become prominent for studying both quasi-static \cite{francfort1998revisiting,bourdin2008variational,bourdin2000numerical,yoon2021quasi,lee2022finite,burke2013adaptive,burke2010adaptive} and dynamic \cite{larsen2010models,larsen2010existence,bourdin2011time} evolution problems. In these phase-field methods, the sharp crack surface is approximated by a smooth, scalar-valued function (with values close to $0$ in damaged regions and $1$ in undamaged areas), which introduces a transition region between the damaged and undamaged zones. The huge advantage of such approaches is that crack initiation, nucleation, evolution, kinking, branching, and the formation of curvilinear paths are automatically embedded in the model. Most notably, the post-processing of local crack criteria, such as calculating stress intensity factors and the need for remeshing the entire domain, is not required.

The foundational theory of brittle fracture proposed by Francfort and Marigo recasts the complex problem of crack propagation into a variational energy minimization framework \cite{francfort1998revisiting}. At its core, the model posits that a crack initiates and evolves to minimize the system's total energy at every moment. This energy is composed of two competing parts: the stored elastic energy within the bulk material and the dissipated energy required to form new crack surfaces. A direct computational treatment of the sharp cracks in this model is notoriously difficult. To make the problem more tractable, the regularized version known as the Ambrosio-Tortorelli phase-field model is often employed. This powerful mathematical technique approximates the sharp discontinuity by introducing a continuous ``phase-field'' variable, which smoothly transitions between the undamaged and fully broken states of the material. Even with this regularization, the direct, monolithic solution of the resulting full coupled, non-linear Euler–Lagrange equations presents significant stability and convergence problems. Therefore, these systems are typically addressed using sophisticated numerical approaches, such as staggered solution schemes, line-search algorithms, or other advanced methods to ensure robust and accurate results \cite{gerasimov2016line,miehe2010phase,amor2009regularized,natarajan2019phase}.

While the length scale parameter is a critical component for regularizing the sharp crack topology in phase–field models (PFMs) of brittle fracture, its introduction raises several significant challenges. A primary issue concerns its physical interpretation; it remains a subject of debate whether it is a purely numerical regularization parameter or a genuine material property representative of the material's microstructure \cite{wu2018length,bourdin2008variational}. In many classical PFMs, the length scale directly governs the predicted tensile strength of the material, creating an undesirable dependence of a physical property on a numerical parameter \cite{miehe2010thermodynamically}. This dependency is intrinsically linked to the computational cost. Accurately resolving the diffuse damage zone necessitates a mesh size significantly smaller than the length scale. Consequently, approximating a sharp crack by choosing a very small length scale leads to prohibitive computational expense, particularly for large-scale, three-dimensional simulations \cite{borden2012phase}. Furthermore, standard PFMs often suffer from spurious damage initiation at stress concentrations, where damage may evolve prematurely before the global failure criterion is met \cite{amor2009regularized}. This artifact, which also depends on the length scale, can compromise the predicted stiffness and peak load of the structure. In response to these limitations, significant research efforts are directed towards developing more advanced formulations—such as those employing hybrid or unified degradation functions—that aim to decouple the regularization length scale from the physical fracture response \cite{wu2017unified}. In recent work \cite{phansalkar2022spatially,phansalkar2023extension}, a model is proposed that extends the classical Ambrosio-Tortorelli (\textsf{AT2}) framework by including the damage length scale parameter as a third unknown variable in the global energy minimization. However, the foundational \textsf{AT2} model possesses a well-documented limitation \cite{pham2011gradient}. Specifically, its formulation is nonlinear in elastic variation prior to fracture, which induces a non-physical artifact of early material degradation, causing the model to predict a loss of stiffness at low stress levels where the material should behave purely elastically.

In this work, we extend previous studies by incorporating the damage length scale parameter as a third primary variable within the classical Ambrosio-Tortorelli (\textsf{AT1}) phase-field framework, in addition to the mechanical displacement and regularized phase-field variables. The total energy functional, composed of three distinct terms, is minimized with respect to all three variables simultaneously. A key innovation of this approach is the utilization of the dynamically optimized and spatially varying length parameter as a natural and physically motivated indicator for guiding local adaptive mesh refinement. To solve the coupled system of partial differential equations derived from the Euler-Lagrange equations of the total energy functional, a conforming Continuous Galerkin finite element method is employed. Our results demonstrate that the optimal length parameter, whether treated as a global or locally varying field, is significantly larger than the values typically used in standard regularized phase-field models. Such a large parameter is doubly advantageous: it reduces computational cost by allowing for a coarser mesh discretization, and it serves as a key factor in obtaining robust numerical convergence for the approximation.

This paper is organized as follows. Section~\ref{intro_model} introduces the governing model for quasi-static brittle fracture and presents a modified energy functional designed to compute an optimal regularization length parameter. We formulate two distinct variants of this model: the first determines a single, optimal uniform length parameter for the entire domain, while the second treats the length scale as a spatially varying field, which is solved for concurrently with the displacement and phase fields. In Section~\ref{fem}, we detail the finite element discretization of the proposed three-field problem, presenting both its continuous variational form and its final discrete algebraic formulation.  Section~\ref{rd} is dedicated to the numerical implementation and results. It includes discussions on key numerical aspects, parameter calibration, and the staggered iterative algorithm employed for the computations. Subsequently, we present and analyze numerical results for benchmark tests involving anti-plane shear loading. Finally, Section~\ref{conclusions} summarizes the key findings of this work and outlines potential directions for future research.

\section{Mathematical Model for Quasi-Static Fracture}\label{intro_model}

Let $\bfa{X}(t_0) \in \mathbb{R}^2$ be an arbitrary point in the stress-free reference configuration ${\Omega}_{t_0}$ of a material body. The corresponding point in the current configuration ${\Omega}_t$ at time $t$ is $\bfa{x}(t) = \bfa{f}(\bfa{X}(t_0), t)$, where $\bfa{f}$ denotes the deformation map. The displacement vector $\bfa{u}(\bfa{x})$ and the deformation gradient $\bfa{F}(\bfa{x})$ are defined by
\begin{equation}
\bfa{u}(\bfa{x}) := \bfa{x} - \bfa{X}, \quad \bfa{F} := \frac{\partial\bfa{f}}{\partial\bfa{X}}.
\end{equation}
The left Cauchy-Green tensor $\bfa{B}$, the right Cauchy-Green tensor $\bfa{C}$, and the Green-Saint Venant strain tensor $\bfa{E}$ are defined as
\begin{equation}
\bfa{B} := \bfa{F}\bfa{F}^\textit{T}, \quad \bfa{C}:= \bfa{F}^\textit{T}\bfa{F}, \quad \bfa{E} := \frac{1}{2}(\bfa{C}-\bfa{I}).
\end{equation}
Assuming the displacement gradients are ``small'' in the sense that
\begin{equation}
\max_{\bfa{x} \in \Omega_{t_0}, \; t \in \mathbb{R}} |\nabla \bfa{u}(\bfa{x})| \ll o(\delta), \quad \delta \ll 1. 
 \label{LAssumption}
\end{equation}
Linearizing $\bfa{E}$ under assumption \eqref{LAssumption}, the nonlinear strain measure $\bfa{E}$ reduces to the linearized strain tensor $\bfeps$:
\begin{equation}
\bfeps(\bfa{x}) = \frac{1}{2} \left( \nabla \bfa{u}(\bfa{x}) + \nabla \bfa{u}^{T}(\bfa{x}) \right).
\end{equation}
Consequently, $\bfa{E}(\bfa{x})$ can be approximated by $\bfeps(\bfa{x})$. 
The problem considered in this study is crack propagation under a quasi-static state of \textbf{anti-plane shear} in a 2D domain. In anti-plane problems, all field quantities (such as the displacement vector $\bfa{u}(\bfa{x})$, stress tensor $\bfsigma(\bfa{x})$, and strain tensor $\bfeps(\bfa{x})$) depend only on the in-plane spatial variables $x_1$ and $x_2$. The displacement vector $\bfa{u}(\bfa{x})$ is characterized by a single non-zero component, $u(x_1, x_2)$, in the out-of-plane direction (herein denoted $x_3$):
\begin{equation}\label{eq:dis}
\bfa{u}(\bfa{x}) = (0, 0, u(x_1, x_2)).
\end{equation}
This specific displacement field signifies that the material undergoes \textbf{pure shear deformation}. The only non-zero components of the infinitesimal strain tensor $\bfeps(\bfa{x})$ are the out-of-plane shear strains $\epsilon_{13} = \epsilon_{31} = \frac{1}{2} \frac{\partial u}{\partial x_1}$ and $\epsilon_{23} = \epsilon_{32} = \frac{1}{2} \frac{\partial u}{\partial x_2}$. All normal strains ($\epsilon_{11}, \epsilon_{22}, \epsilon_{33}$) and the in-plane shear strain ($\epsilon_{12}$) are identically zero. Consequently, for a linear isotropic elastic solid, the material experiences a state of \textbf{pure shear stress}. The only non-zero components of the stress tensor $\bfsigma(\bfa{x})$ are the corresponding out-of-plane shear stresses $\sigma_{13} (=\sigma_{31})$ and $\sigma_{23} (=\sigma_{32})$. The constitutive relationship between the stress and strain tensors simplifies to
\begin{equation}\label{Eq:Constitutive}
\bfsigma = 2 \mu \bfeps,
\end{equation}
where $\mu$ is the shear modulus. This tensorial equation directly yields the component-wise relations for the non-zero shear stresses: $\sigma_{13} = 2\mu\epsilon_{13} = \mu \frac{\partial u}{\partial x_1}$ and $\sigma_{23} = 2\mu\epsilon_{23} = \mu \frac{\partial u}{\partial x_2}$.

\subsection{Variational framework for quasi-static crack evolution } 
This study investigates crack propagation within a deformable body $\Omega$, an open, bounded, and connected subset of $\mathbb{R}^{2}$ featuring a Lipschitz boundary $\partial \Omega$. Our analysis employs standard notations for Lebesgue spaces, denoted $L^{p}(\Omega)$ for $p \in [1, \infty)$, and Sobolev spaces like $H^1(\Omega)$.

Central to modeling discontinuities, such as cracks, is the concept of functions with \textbf{bounded variation} ($BV$). A function $f \in L^1(\Omega)$ belongs to $BV(\Omega)$ if its total variation is finite. This condition is met if the following supremum is bounded \cite{burke2010adaptive, burke2013adaptive}:
\begin{equation}
\sup \left\{ \int_{\Omega} f \, \operatorname{div} \boldsymbol{\varphi} \, d\bfa{x} : \boldsymbol{\varphi} \in C_c^1(\Omega; \mathbb{R}^{N}), \|\boldsymbol{\varphi}\|_{L^\infty} \leq 1 \right\} < \infty.
\end{equation}
Functions in $BV(\Omega)$ are notable for admitting jump discontinuities, which are captured by their distributional gradients.

A refined space, crucial for fracture mechanics, is that of \textbf{special functions of bounded variation}, $SBV(\Omega)$, a subset of $BV(\Omega)$ \cite{ambrosio1990metric}. The distributional derivative $Df$ for any $f \in SBV(\Omega)$ decomposes into two parts: an absolutely continuous part and a jump part:
\begin{equation}
Df = \nabla f \mathcal{L}^N + (f^+(x) - f^-(x)) \otimes \boldsymbol{\nu}_f(x) \mathcal{H}^{N-1} \lfloor S(f).
\end{equation}
Here, $\nabla f$ represents the approximate gradient of $f$, while $S(f)$ denotes its jump set (i.e., the set where $f$ is not a Lebesgue point). The vector $\boldsymbol{\nu}_f$ is the unit normal to $S(f)$, and $f^+$ and $f^-$ are the one-sided trace values of $f$ at $S(f)$. The $N$-dimensional Lebesgue measure is $\mathcal{L}^N$, and $\mathcal{H}^{N-1}$ is the $(N-1)$-dimensional Hausdorff measure.

We consider $\Omega$ as the reference configuration of a linearly elastic material, initially intact and free of stress. A \textbf{crack set}, $\Gamma$, represents an internal discontinuity within $\Omega$. This set $\Gamma$ is understood as a lower-dimensional manifold (a $1$D curve for a two-dimensional curve). It is not necessarily connected but is wholly contained in $\Omega$ and is assumed to be an $\mathcal{H}^{N-1}$-measurable set. Critically, $\Gamma$ corresponds to the jump set $S(\bfa{u})$ of the displacement field $\bfa{u}$.\footnote{The \textit{crack set} $\Gamma$ corresponds to the \textit{jump set} $S(\bfa{u})$ across which the displacement field $\bfa{u}$ is discontinuous.} For problems like anti-plane shear (cf. Eq.~\eqref{eq:dis}), the relevant scalar displacement component $u(\bfa{x})$ is a real-valued function belonging to $SBV(\Omega)$.

The mechanical response of the fracturing material is governed by an energy functional comprising two main contributions:
\begin{itemize}
    \item The \textbf{elastic (or bulk) energy}, $E_b: H^1 \to \mathbb{R} \cup \{ + \infty \}$, stored in the uncracked portion of the body:
    \begin{equation}\label{eq:elastic_revised}
    E_b(\bfa{u}) := \int_{\Omega \setminus S(\bfa{u})} W(\nabla u) \, d\bfa{x}.
    \end{equation}
    For the specific case of anti-plane shear, where the displacement vector $\bfa{u}=(0,0,u)$, $W(\nabla u)$ simplifies to $\|\nabla u\|^2$, with $u$ being the scalar out-of-plane displacement and $\nabla$ its in-plane gradient. More generally, $W(\boldsymbol{\varepsilon}(\bfa{u}))$ is the strain energy density.
    \item The \textbf{surface energy}, $E_s$, associated with the creation or presence of cracks:
    \begin{equation}\label{eq:surene_revised}
    E_s(\Gamma) := \int_{\Gamma} \kappa(\bfa{x}) \, d\mathcal{H}^{N-1}(\bfa{x}).
    \end{equation}
    Here, $\kappa(\bfa{x}) > 0$ denotes the material's \textbf{fracture toughness} or critical energy release rate (often symbolized as $G_c$) \cite{francfort1998revisiting, bourdin2000numerical, bourdin2008variational, yoon2021quasi,manohar2025convergence}. A value of $\kappa(\bfa{x}) = +\infty$ implies the material is unbreakable at $\bfa{x}$. For the current work, $\kappa$ is often taken as a positive constant. If so, $E_s(\Gamma) = \kappa \mathcal{H}^{N-1}(\Gamma)$.
\end{itemize}
The \textbf{total energy} of the system $(\bfa{u}, \Gamma)$ is then postulated as:
\begin{equation}\label{eq:totene_revised}
E(\bfa{u}, \Gamma) :=
\begin{cases}
E_b(\bfa{u}) + E_s(\Gamma) & \text{if } \mathcal{H}^{N-1}(S(\bfa{u}) \setminus \Gamma) = 0 \\
+\infty & \text{otherwise.}
\end{cases}
\end{equation}
The constraint $\mathcal{H}^{N-1}(S(\bfa{u}) \setminus \Gamma) = 0$ ensures that displacement discontinuities $S(\bfa{u})$ are confined to the designated crack path $\Gamma$. The model accommodates crack growth, typically governed by an irreversibility constraint stating that cracks, once formed, do not heal.\footnote{Irreversibility implies that once a crack forms or extends, it cannot heal spontaneously.}

Quasi-static crack evolution is driven by time-varying external loads, implemented through Dirichlet boundary conditions $g(t)$. These conditions are prescribed on an open portion $\partial\Omega_D$ of the boundary $\partial\Omega$, while the remainder, $\partial\Omega_N = \partial\Omega \setminus \overline{\partial\Omega_D}$, is traction-free. The function $g(t)$ is assumed to satisfy $g \in L^{\infty}(0, T; W^{1, \infty}(\Omega)) \cap W^{1,1}(0, T; H^1(\Omega))$. For each $g(t)$, the set of kinematically admissible displacements is
\begin{equation}\label{eq:bc_revised}
\mathcal{A}(g(t)) := \{ \bfa{u} \in SBV(\Omega; \mathbb{R}^3) \text{ such that } u_3 \in SBV(\Omega) : \bfa{u}|_{\partial\Omega_D} = \mathbf{g}(t)|_{\partial\Omega_D} \},
\end{equation}
where $\mathbf{g}(t)=(0,0,g(t))$ for anti-plane shear scenarios. The core of the quasi-static fracture model involves finding configurations $(\bfa{u}, \Gamma)$ that minimize the total energy $E$ at discrete load steps (or pseudo-time instances) $t$, subject to the boundary conditions and the crack irreversibility principle \cite{francfort2003existence}. This variational approach was pioneered by Francfort and Marigo \cite{francfort1998revisiting} and subsequently advanced, notably by Bourdin \cite{bourdin2000numerical}.

The loading history over an interval $[0, T]$ is discretized into $0 = t_0 < t_1 < \dots < t_M = T$, where $M$ is the number of load increments, and $\Delta t_{\max} = \max_{k=1,\dots,M} (t_k - t_{k-1})$ is the largest increment. Starting from a known initial state $(\bfa{u}(t_0), \Gamma(t_0))$, where $\Gamma(t_0)$ is the jump set $S(\bfa{u}(t_0))$, the problem is to determine the evolving displacement-crack pair $(\bfa{u}(t_j), \Gamma(t_j))$ for each subsequent load step $t_j$, where $j=1, \dots, M$.

\subsection{Energy Minimization for the Francfort-Marigo Model}

This section provides a brief overview of the minimization algorithm used to compute the quasi-static fracture evolution according to the Francfort-Marigo model. The solution algorithm is formulated at the continuous level, followed by a brief overview of its finite element implementation. The algorithm presented is a standard approach and has been widely adopted in the literature \cite{yoon2021quasi,bourdin2008variational,burke2010adaptive,lee2022finite,heister2015primal,lee2016pressure,amor2009regularized}.

For each time step $t_j$, where $j = 1, \dots, N$, the displacement $u(t_j)$ and crack set $\Gamma(t_j)$ are determined by the following two-step iterative procedure:
\begin{align}
u(t_j) &:= \argmin_{\mathbf{v} \in \mathcal{A}(g(t_j))} \left\{ E_{b}(\mathbf{v}) + E_{s}(S(\mathbf{v}) \cup \Gamma(t_{j-1})) \right\}, \\
\Gamma(t_j) &:= S(u(t_j)) \cup \Gamma(t_{j-1}).
\end{align}
This procedure enforces the irreversibility of crack propagation, as the crack set can only grow or remain the same. The direct numerical solution of this discrete-time minimization problem is challenging for standard numerical techniques, such as the finite element or boundary element methods. This difficulty arises from the need to explicitly determine the evolving crack set $\Gamma$ and identify the associated displacement jump set $S(u)$. A common and practical approach to overcome this challenge is to employ a regularization of the total energy \eqref{eq:totene_revised}, leading to a formulation amenable to standard numerical implementation.
To this end, a widely adopted regularization is the Ambrosio-Tortorelli functional \cite{ambrosio1990approximation}. This regularized energy $E_{\xi} : \Hspace{1}(\Omega; \mathbb{R}^2) \times \Hspace{1}(\Omega; [0, 1]) \to \mathbb{R}$ (where $d$ is the spatial dimension) is defined for parameters $\eta > 0$ and $\xi > 0$ as follows:
\begin{equation}\label{reg:energy}
E_{\xi} (u, v) := \int_{\Omega} \left( (1 - \eta) v^2 + \eta \right) \|\nabla u\|^2 \, d \mathbf{x} + \frac{G_c}{c_v} \int_{\Omega} \left[ \frac{(1-v)}{\xi} + \xi \|\nabla v \|^2 \right] \, d\mathbf{x}.
\end{equation}
In \eqref{reg:energy}, the phase-field variable $v \in \Hspace{1}(\Omega; [0, 1])$ is a smooth scalar function. It approximates the crack set $\Gamma$: $v \approx 1$ in the undamaged material and $v \approx 0$ on the crack. The parameter $\xi > 0$ represents a regularization length, controlling the width of the diffuse interface approximating the crack. The symbol $\| \cdot \|$ is the classical $L^2$-norm. 
The term $ \left[ \left( (1 - \eta) v^2 + \eta \right)  \| \nabla u \|^2 \right]$ in \eqref{reg:energy} penalizes displacement gradients less severely where $v \approx 0$ (i.e., on the crack), effectively modeling the loss of material stiffness. The terms involving $\xi$ approximate the surface energy and ensure that $v$ transitions sharply from $1$ to $0$ over a width proportional to $\xi$.
The minimization of $E_{\xi}$ with respect to both $u$ and $v$ provides an approximation consistent with Griffith's criterion for fracture as $\xi \to 0$ (in the sense of $\Gamma$-convergence).

The sequence of regularized functionals $\{E_{\xi}\}_{\xi > 0}$, defined by \eqref{reg:energy}, is known to $\Gamma$-converge to the original Griffith energy functional $E$ (given in \eqref{eq:totene_revised}) in the $\Lspace{1}(\Omega) \times \Lspace{1}(\Omega)$ topology as $\xi \to 0$ \cite{braides2002gamma}. For any fixed $\xi > 0$ and $\eta > 0$, the existence of minimizers for $E_{\xi}$ was established in \cite{AmTo92}. The regularization parameter $\eta > 0$ is introduced primarily to ensure the well-posedness of the elastic energy term, particularly in regions where the phase field $\varphi$ approaches zero (i.e., on the crack). It is typically chosen as a small, constant positive value throughout the computation. This helps prevent an overestimation of the bulk (strain) energy, which could otherwise lead to an underestimation of the crack surface energy.
While some analyses, such as \cite{bourdin2000numerical} , coupled $\eta$ with $\xi$ (e.g., $\eta \propto \xi$) to maintain the validity of $\Gamma$-convergence for specific problem settings, the original Ambrosio-Tortorelli $\Gamma$-convergence framework effectively targets a limiting functional where the elastic energy in the fully cracked state (where $\varphi=0$) is zero, consistent with an unregularized behavior in the limit.
For quasi-static simulations, the optimal choice of $\eta$ that yields physically meaningful crack patterns consistent with experimental observations can be problem-dependent and remains an area of ongoing research.

The regularized quasi-static minimization problem is then formulated as follows:
Given initial conditions $u_{\xi}(t_0)$ and $v_{\xi}(t_0) = 1$ almost everywhere (representing an initially intact material), find the displacement $u_{\xi}(t_k)$ and phase field $v_{\xi}(t_k)$ for each discrete time $t_k$, for $k=1, \dots, N$ (where $N$ is the total number of time steps), by solving the following minimization problem at each step:
\begin{align}\label{min:ene:formulation}
(u_{\xi}(t_k), v_{\xi}(t_k)) = \argmin \bigg\{ E_{\xi}(\widehat{u}, \widehat{v}) \,:\,
& \widehat{u} \in \Hspace{1}(\Omega; \mathbb{R}^2), \widehat{u}|_{\partial\Omega_D} = \mathbf{g}(t_k); \nonumber \\
& \widehat{v} \in \Hspace{1}(\Omega; [0,1]), \widehat{v}(\mathbf{x}) \leq v_{\xi}(\mathbf{x}, t_{k-1}) \text{ for a.e. } \mathbf{x} \in \Omega
\bigg\}.
\end{align}
The constraint $\widehat{v}(\mathbf{x}) \leq v_{\xi}(\mathbf{x}, t_{k-1})$ enforces the irreversibility of damage, meaning the material cannot heal.

This formulation for quasi-static evolution implies that the crack (or damaged region) is approximated by a narrow zone where the phase field $v_{\xi}(t_k)$ is close to zero. Conversely, the undamaged portion of the material corresponds to regions where $v_{\xi}(t_k)$ is close to one.
The width of the transition layer between these states is governed by the regularization parameter $\xi$. As $\xi \to 0$, the phase-field approximation is expected to sharpen, ideally converging to a discrete crack representation where $v = 1$ almost everywhere in the undamaged domain and $v = 0$ on the lower-dimensional crack set.

A key challenge in the aforementioned quasi-static phase-field formulations is the selection of the regularization length parameter, $\xi$. Its optimal determination often lacks a definitive consensus in the existing literature. Addressing this, Phansalkar et al.~\cite{phansalkar2022spatially, phansalkar2023extension} introduced adaptive models where $\xi$ is treated as an additional field variable, and the total energy functional is minimized with respect to displacement, the phase field, and this spatially varying regularization length.

The present study extends these adaptive concepts by incorporating a linear material degradation function (characteristic of \textsf{AT1}-type models) and specifically considers its application under anti-plane shear loading conditions.
A critical aspect of such adaptive models is ensuring the boundedness of the regularization length field, $\xi$. Therefore, this work proposes a revised energy functional that incorporates two positive parameters designed to maintain $\xi$ within desired upper and lower bounds. The modified energy functional is consequently formulated as:
\begin{align}
E(u(x), \,  v(x), \, \xi) = &
\frac{\mu}{2} \int_{\Omega} \Big((1 - \eta)v^2+\eta\Big) \| \nabla u \|^2  \;  d \mathbf{x}  + \frac{G_c }{c_v} \int_{\Omega} \left( \frac{1-v+\zeta}{\xi} + \xi \| \nabla v(x) \|^2 \right)  \; d \mathbf{x} \notag \\
 & + \int_{\Omega} \alpha \; \xi \, d\bfa{x}. 
\end{align}
In the proposed energy formulation, the parameters $\zeta$ and $\alpha$ are implicitly dependent on $\xi$. This formulation aims to employ a smaller $\xi$ value in the vicinity of the crack-tip, while utilizing a marginally larger value in other regions. Such an approach offers an advantage over methods that use a globally constant $\xi$, which typically require $\xi > h$ within the damage zone. The stability and convergence of solutions to the Francfort-Marigo model are known to be critically dependent on the $\xi$ parameter. Hence, choosing an optimal value of $\xi$ is crucial for obtaining stable solutions to the above model.  The proposed energy is minimized with respect to the three field variables: displacement $u(x)$, phase-field $v(x)$, and length scale parameter $\xi$. The minimization concerning $\xi$ can be performed assuming either a globally constant value or a spatially varying field.

Before discussing the minimization of the total energy functional, it's essential to define the appropriate function spaces at a specific time instant $t=t_j$. We consider $\mathcal{V}=H^1(\mathcal{D})$ as the fundamental Sobolev space. From this, we define the pertinent subspaces as:
\begin{subequations}
\begin{align}
& \mathcal{V}_d:=\big\{\varphi \in \mathcal{V}\,|~~\varphi=0~~ \text{on}~~ \Omega_d \big\},  \label{eq:Vd_def_rephrased} \\
& \mathcal{V}_c:=\big\{\varphi \in \mathcal{V}\,|~~ \varphi=0 ~~ \text{on}~~ CR(t_{j-1}) \big\}, \label{eq:Vc_def_rephrased} \\
\text{and} ~~~~ &\nonumber\\
& \mathcal{V}_f:=\big\{\varphi \in \mathcal{V}\,|~~ \varphi=f(t_j)~~ \text{on}~~ \Omega_d \big\}, \label{eq:Vf_def_rephrased}
\end{align}
\end{subequations}
In these expressions, $\Omega_d$ signifies the Dirichlet boundary, $f(t_j)$ is a given function on this boundary at time $t_j$, and $CR(t_{j-1})$ represents the crack set (i.e., the irreversibly damaged zone) determined at the preceding time step $t_{j-1}$. The crack set for the current time $t_j$ is defined as
\begin{equation}
CR(t_j) = \left\{ x \in \overline{\mathcal{D}} \; \colon \; v_{\epsilon} (x, \; t_j) < TOL \right\}, \label{eq:CR_def_rephrased}
\end{equation}
where $TOL$ denotes a predefined small positive tolerance value. To enforce crack irreversibility, the phase-field variable $v_{\epsilon}(x, t_j)$ is set to zero for any point $x$ located within this crack set $CR(t_j)$ (which is assumed to be non-empty).

\textbf{A constant, globally optimal value for $\xi$: }
First, the analysis considers a spatially uniform optimal length scale parameter, $\xi$, within the domain $\Omega$. Minimization of the energy functional with respect to the displacement field $u(x)$, the phase-field $v(x)$, and the constant parameter $\xi$ yields the following system of Euler-Lagrange equations:
\begin{subequations}
\begin{align}
    E_u &= \mu\int_{\Omega} \left( (1-\eta) v ^2 + \eta \right) \nabla u \cdot \nabla {\phi} \, dx  &&  \text{in } \Omega \label{eq:Eu_const_xi} \\[10pt]
    E_v &= \mu \int_{\Omega} \| \nabla u \|^2 (1-\eta) v \phi \, dx 
         - \frac{G_c}{c_v} \int_{\Omega} \frac{\phi}{\xi}\, dx 
         + \frac{2G_c}{c_v} \int_{\Omega} \xi \nabla v \cdot \nabla \phi \, dx  &&  \text{in } \Omega \label{eq:Ev_const_xi} \\[10pt]
    E_{\xi} &= \frac{G_c}{c_v} \int_{\Omega}\left( \frac{-(1 - v+\zeta)}{\xi^2} + \| \nabla v \|^2 \right) \, dx 
            + \int_{\Omega} \alpha \, dx &&  \text{in } \Omega \label{eq:Exi_const_xi} \\[10pt]
    \xi &= \sqrt{ \frac{ \frac{G_c}{c_v} \int_{\Omega} (1-v+\zeta)\, dx}
                 {\frac{G_c}{c_v}\int_{\Omega} \| \nabla v \|^2 dx + \int_{\Omega} \alpha \,dx }} \label{eq:xi_const_xi}
\end{align}
\end{subequations}
In these equations, $E_u$, $E_v$, and $E_{\xi}$ represent the variations of the energy functional with respect to $u$, $v$, and $\xi$, respectively. The parameter $\mu$ is an elastic modulus, $k$ is a residual stiffness parameter, $G_c$ is the critical energy release rate, $c_v$ is a model constant, and $\eta$ is a regularization parameter. The variable $\phi$ denotes a generic test function (with $\boldsymbol{\phi}$ for its vector counterpart if $u$ is a vector field). The optimal constant $\xi$ is determined from $E_{\xi}=0$.

\textbf{A spatially dependent (or an inhomogeneous) $\xi$: }Subsequently, the model is extended by treating $\xi$ as a spatially varying field variable, denoted $\xi(x)$. Minimization of the modified energy functional with respect to $u(x)$, $v(x)$, and $\xi(x)$ leads to the weak-form Euler-Lagrange equations for $u(x)$ and $v(x)$, and a pointwise algebraic optimality condition for $\xi(x)$:
\begin{subequations}
\begin{align}
    E_u &= \mu \int_{\Omega} \left( (1-\eta) v^2 + \eta \right) \nabla u \cdot \nabla {\phi} \, dx     &&  \text{in } \Omega \label{eq:Eu_field_xi} \\[10pt]
    E_v &= \mu \int_{\Omega} \| \nabla u \|^2 (1-\eta) v \phi \, dx 
         - \frac{G_c}{c_v} \int_{\Omega} \frac{\phi}{\xi(x)}\, dx 
         + \frac{2G_c}{c_v} \int_{\Omega} \xi(x)  \nabla v \cdot \nabla \phi \, dx   &&  \text{in } \Omega \label{eq:Ev_field_xi} \\[10pt]
    E_{\xi} &= \frac{G_c}{c_v} \left( \frac{-(1 - v(x)+\zeta)}{\xi(x)^2} + \| \nabla v(x) \|^2 \right) 
             +  \alpha &&  \text{pointwise in } \Omega \label{eq:Exi_field_xi} \\[10pt]
    \xi(x) &= \sqrt{ \frac{ \frac{G_c}{c_v} (1-v(x)+\zeta)} 
                 {\frac{G_c}{c_v} \| \nabla v(x) \|^2 + \alpha }}  && \text{pointwise in } \Omega \label{eq:xi_field_xi}
\end{align}
\end{subequations}
Here, the equations for $E_u$ and $E_v$ remain integral forms. The condition $E_{\xi}=0$ (Eq.~\eqref{eq:Exi_field_xi}) provides a local algebraic equation for $\xi(x)$, whose solution is given by Eq.~\eqref{eq:xi_field_xi}, assuming $\alpha$ can also be a spatially varying function $\alpha$.

\section{Finite-element discretization}\label{fem}
We begin by outlining the discrete finite element framework, which utilizes a conforming simplicial mesh, denoted by $\Th$, to partition the domain closure $\overline{\Omega}$.
This mesh accurately represents the domain boundary $\partial \Omega$.
Furthermore, for any two distinct simplices $\tau_i, \tau_j \in \Th$ (with $i \neq j$), their intersection is restricted to being either the empty set, a common vertex, or a shared edge.
The characteristic diameter of each simplex $\tau \in \Th$ is denoted by $h_\tau = \text{diam}(\tau)$, and the overall mesh resolution is given by $h = \max_{\tau \in \Th} h_\tau$.
The mesh family $\Th$ is assumed to satisfy a shape-regularity condition; this implies the existence of a positive constant $c_{\varrho}$ such that for every simplex $\tau \in \Th$, the inequality $h_\tau / \varrho_\tau \leq c_{\varrho}$ holds, where $\varrho_\tau$ is the diameter of the largest circle that can be inscribed within $\tau$.
Each simplex $\tau \in \Th$ is also affinely equivalent to a reference simplex $\refface$, meaning that for each $\tau$, there exists an invertible affine map $M_{\tau}: \refface \to \tau$.The reference simplex $\refface$ is defined as $$\refface = \{ \refcoord \in \mathbb{R}^2 : \hat{z}_k > 0 \text{ for } k=1,\dots,d, \text{ and } \sum_{k=1}^d \hat{z}_k < 1 \}.$$
Lastly, corresponding to the vertices $\{ \nodecoord_i \}_{i \in \nodes}$ of the mesh $\Th$, where $\nodes$ is the index set for these vertices, we employ a nodal basis.
This basis is comprised of functions $\basis_i$ for $i \in \nodes$ that are continuous across $\overline{\Omega}$, piecewise linear with respect to $\Th$, and satisfy the Lagrange property $\basis_i(\nodecoord_j) = \delta_{ij}$ (the Kronecker delta) for all $j \in \nodes$.

Based on this framework, we define the following finite-dimensional approximation spaces at a discrete time $t=t_n$:
\begin{align*}
\FESpace &:= \text{span} \{ \basis_i : i \in \nodes \} = \left\{ \phi_h = \sum_{i \in \nodes} \Phi_i \basis_i : \Phi_i \in \mathbb{R} \right\}, \\
\FESpaceHomD &:= \left\{ \phi_h \in \FESpace : \Phi_i = 0 \text{ for all } i \in \nodesDboundary \right\}, \\
\FESpaceInhomD &:= \left\{ \phi_h \in \FESpace : \Phi_i = g(t_n, \nodecoord_i) \text{ for all } i \in \nodesDboundary \right\},
\end{align*}
where $\nodesDboundary = \{i \in \nodes : \nodecoord_i \in \partial \Omega \}$ is the index set of nodes on the Dirichlet boundary, and $g(t_n, \nodecoord_i)$ represents the prescribed Dirichlet data.

To handle evolving crack patterns, we define a discrete crack set based on a phase-field variable $\phasefield \in \FESpace$ from the preceding time step $t_{n-1}$ and a tolerance parameter $\tolCR$:
The set of faces considered "cracked" at $t_{n-1}$ is
$$ \FacesCR(t_{n-1}) := \left\{ f \in \Fh : \phasefield(\nodecoord, t_{n-1}) \leq \tolCR \text{ for all } \nodecoord \in \bar{f} \right\}. $$
The discrete crack region is then the union of these faces: $\CrackRegion(t_{n-1}) := \bigcup_{f \in \FacesCR(t_{n-1})} \bar{f}$.
This allows the definition of a finite element space that incorporates crack irreversibility from the previous step:
$$ \FESpaceCrack := \left\{ \phi_h \in \FESpace : \phi_h(\nodecoord) = 0 \text{ for all } \nodecoord \in \CrackRegion(t_{n-1}) \right\}. $$
For brevity, when the specific time $t_n$ is understood, the notations $\FESpaceCrack$ and $\FESpaceInhomD$ may be simplified to $\FESpaceCrackSimple$ and $\FESpaceInhomDSimple$, respectively.

\subsection{Continuous formulation}
The continuous weak formulation for quasi-static crack evolution using the proposed three-field approach is obtained by taking the Lagrangian of the total energy functional. Such an operation yields the $L^1$-formulation of the total energy system. This system arises from the stationarity conditions of an energy functional that includes both the elastic energy of the deforming body and the energy associated with crack formation, regularized via a phase-field variable. This phase-field variable smoothly transitions between an intact and a fully broken material state. Such a formulation is posed in $\mathcal{V}$ seeking the solution in terms of three-field variables: $u, \; v$, and $\xi$ (both optimal uniform global length variable and spatially-varying length variable). A key aspect is often the enforcement of crack irreversibility, meaning the phase-field variable representing damage cannot heal.

For $n =1, \, 2, \, \ldots$, the problem seeks an approximate solution $u_n \in \mathcal{V}_f$ and $v_n \in \mathcal{V}_c$such that:
\begin{equation}
    a_i(\phi_h, \varphi_h) = l_i(\varphi_h),  \quad i=1,2, \\
\end{equation}
where $a_i(\cdot,\cdot)$ is the bilinear form, and $l(\cdot)$ represents the linear functional associated with the weak formulation. The particular definitions of these terms are as follows:
\begin{subequations}
\begin{align}
 a_1(u_n, \varphi) &= \int_\Omega\mu \left( (1-\eta) v^2 + \eta \right) \nabla u \cdot \nabla \phi d\bfa{x}, \quad \forall \; \phi \in \mathcal{V}_d \\
 l_1(\varphi)  &=0 \\
  a_2(v_n, \varphi) &=\int_\Omega \mu  ||\nabla u||^2 (1-\eta) v \phi + \frac{2G_c}{c_v} \xi \nabla v \cdot \nabla \phi d\bfa{x}, \quad \forall \; \phi \in \mathcal{V}\\
  l_2(\varphi) &=\frac{G_c}{c_v}\int_\Omega  \frac{\phi}{\xi} d\bfa{x}, \quad \forall \; \phi \in \mathcal{V}
\end{align}
\end{subequations}
The integral formulation presented above establishes the well-posedness of these equations.
However, since exact solutions are generally unavailable, numerical methods are required.
A collocation approximation, such as the finite element method, is a preferred choice for obtaining these numerical solutions. The subsequent subsection will detail the discrete formulation derived from this continuous weak formulation.

\subsection{Discrete formulation}
This section introduces a discrete formulation for simulating the quasi-static evolution of cracks within a linear elastic material.
The primary objective is to calculate the solution fields, $u^n$ and $v^n$, at a specific time instant $t=t_n$.
Subsequently, the methodology aims to determine an optimal global damage characteristic length and a locally varying characteristic length.

\begin{dwf}
Let $\mu, \; \eta, \; \zeta, \; \alpha$ and a mesh size $h$ are all known and given the boundary conditions for $u(t_n) = g(t_n)$ on $\partial \Omega_D$, then for $n=1, \, 2, \, \ldots$ find $u^n \in \FESpaceInhomD$ and $v^n \in \FESpace$ such that 
\begin{subequations}
\begin{align}
\sum_{K \in T_h} \int_K \mu\left((1 - \eta)v_h^2 + \eta\right) \nabla u^n_h \cdot \nabla \varphi_h \,d\bfa{x}  &= 0, \quad \forall \varphi_h \in \FESpaceHomD  \\
\sum_{K \in T_h} \int_K \mu ||\nabla u_h||^2 (1 - \eta)v_h  \varphi_h \,d\bfa{x}  &+ \frac{2 G_c}{c_v} \sum_{K \in T_h} \int_K \xi \nabla v_h \cdot \nabla \varphi_h \,d \bfa{x}  \notag \\
&=\frac{G_c}{c_v} \sum_{K \in T_h} \int_K \frac{\varphi_h}{\xi} \,d\bfa{x}, \quad \forall \varphi_h \in \FESpace
\end{align}
\end{subequations}
\end{dwf}
Once the discrete weak formulations are assembled, the resulting matrix system is addressed using a direct solver. Following this, the optimal characteristic length for the damage zone is calculated based on the formulas provided in the preceding subsection.

\section{Discussion and numerical considerations} \label{rd}
This section illustrates the theoretical findings through a specific example that employs the adaptive algorithms detailed in the previous section.  The chosen example involves an elastic unit square containing a single crack along one edge, subjected to antiplane shear loading conditions on its boundary.  The choice of an antiplane shear model for a quasi-static crack provides a crucial testbed for our proposed framework. In this configuration, the displacement vector has only one non-zero component, which is perpendicular to the plane of analysis. This essentially enables us to test our formulation in the simplified context of a scalar-valued displacement unknown, isolating the core mathematical and numerical aspects of the problem from the complexities of full vectorial elasticity. We have studied the three-field formulation using a phase-field model with \textsf{AT1}. Instead of modeling a crack as a sharp discontinuity, the \textsf{AT1} model uses a continuous scalar field (the phase-field) that smoothly transitions between an undamaged state (e.g., value of 0) and a fully broken state (e.g., value of 1) over a small length scale. A significant advantage of the \textsf{AT1}  model over its standard counterpart, the \textsf{AT2}  model, is that it predicts a well-defined elastic limit. This means that damage only begins to accumulate after the strain energy reaches a specific, non-zero threshold, which often better represents the behavior of brittle materials.

To solve the coupled system involving the $u$ and $v$ fields, we developed and implemented a multi-step, staggered iterative solver. The quasi-static analysis proceeds by applying a constant displacement increment at each load step, ensuring stable convergence. The irreversibility of the crack, which dictates that a crack cannot heal once formed, is a critical physical constraint that must be enforced numerically. This is typically handled in the code using the following standard technique: The condition is implemented as a local inequality constraint on the phase-field variable, $v$. For any given time step $t_n$, the phase-field value $v(t_n)$ is enforced to be greater than or equal to its value at the previous step, $v(t_{n-1})$. This prevents any decrease in the damage value, effectively making the fracture process irreversible. Within the staggered iterative solver, after computing a potential new phase-field solution at the current step, a simple projection is applied. If any point in the computed field has a lower damage value than in the previous converged step, its value is reset to the previous, higher value before proceeding to the next iteration. Further, if the computed damage value is less than the prescribed small value, then it is set to zero in the subsequent time steps. This could create a crack-widening effect; however, we could not see such an effect in our simulations. To enhance computational efficiency and solution accuracy, this solver is combined with a sophisticated adaptive mesh refinement (AMR) strategy. This AMR scheme employs a multi-step procedure for both \textit{marking} elements that require higher resolution and \textit{coarsening} the mesh in regions of lesser interest.  The entire computational framework was written in \textsf{C++} and built upon the versatile, open-source finite element library \textsf{deal.II} library \cite{arndt2021deal}. 

\subsection{Paramter estimation}
The proposed formulation involves several parameters whose careful selection is critical to ensure the convergence, accuracy, and efficiency of our method. The physical material properties for the simulations were defined by the shear modulus, $\mu=80.8$, and the critical energy release rate, $G_c =2.7$.  Beyond these physical constants, several model-specific parameters must be calibrated. The parameter $\eta$, which regularizes the estimation of the bulk energy, was set to a very small value of $10^{-10}$.  This choice is deliberate, as it prevents the overestimation of surface energy density in the high-gradient region near the crack tip, thereby ensuring the model's physical accuracy.

The model parameters, $\zeta$ and $\alpha$, were calibrated by analyzing the characteristic profile of the phase-field variable, $v$, under different conditions. First, in regions distant from the crack where the material is considered fully intact ($v=1$), and near domain boundaries where the homogeneous Neumann condition ($\mathbf{n} \cdot \nabla v = 0$) causes the phase-field gradient $\| \nabla v \|$ to diminish, a fine resolution of the phase-field is unnecessary. In these areas, we assign a large regularization length scale, specifically $\xi = 10h$, tying it to the local mesh size $h$. This choice reflects a coarsening of the phase-field representation and, when substituted into the model's governing relations, yields a direct expression for the parameter $\zeta$ as a function of the mesh size:
\begin{equation}
    \zeta=\frac{100h^2c_v\alpha}{G_c} \label{eq:eta_calib}
\end{equation}
In contrast, a different approach is required at the crack interface. At the crack tip, the material is fully damaged (i.e., $v=0$), and a sharp transition zone exists where the phase-field rapidly varies. To accurately resolve this high-gradient region, a smaller regularization length that scales with the mesh resolution is necessary. Based on the observation of an inverse proportionality between the gradient norm $\| \nabla v \|$ and the mesh size $h$, we define a smaller regularization length, $\xi = 2h$. Substituting this into the governing equations provides a formula for the parameter $\alpha$:
\begin{equation}
    \alpha=\frac{3G_c}{96c_v h^2} \label{eq:alpha_calib}
\end{equation}
Together, these relationships allow for the systematic calculation of the required parameters for any given mesh discretization. The resulting values for various mesh sizes are summarized in Table \ref{tab:param_calib}.

\begin{table}[h!]
    \centering
    \caption{Calibrated model parameters for different mesh sizes.}
    \label{tab:param_calib}
    \begin{tabular}{cccc}
        \toprule
        {n} & {h} & {${\alpha}$} & {${\zeta}$} \\
        \midrule
        128 & 0.008 & 493.75 & 9.36 \\
        256 & 0.004 & 1975   & 9.36 \\
        512 & 0.002 & 7900   & 9.36 \\
        \bottomrule
    \end{tabular}
\end{table}

\subsection{Staggered algorithm and crack irreversibility}
To solve the coupled system governing phase-field evolution and linear elastic deformation under antiplane shear, a robust staggered, or segregated, solution algorithm is implemented. This computational strategy decouples the full problem into two distinct physics sub-problems that are solved sequentially within each time step or load increment. First, the momentum balance equation is solved for the out-of-plane displacement field, holding the phase-field variable fixed from the previous iteration. In this step, the phase field serves as a static degradation function, locally reducing the material's shear stiffness and effectively representing the current state of material damage. Subsequently, using the newly computed displacement field, the strain energy density is determined and serves as the primary driving force in the second sub-problem: solving the phase-field evolution equation. This update determines the new topology of the diffuse crack network. This alternating solution procedure is repeated until convergence is achieved, providing a computationally efficient and memory-friendly alternative to a monolithic scheme by avoiding the assembly and inversion of a large, fully coupled Jacobian matrix.  This work employs an adaptive iterative scheme to enhance the efficiency of phase-field fracture simulations. The algorithm begins with a coarse mesh and sequentially solves for the displacement, phase-field, and a variable length scale. To ensure accuracy while minimizing computational cost, the mesh is dynamically refined in regions with high solution gradients, enabling a robust and efficient simulation of crack propagation. The overall staggered algorithm is presented in Algorithm~$1$. 

\begin{algorithm}[H]
\caption{Staggered solution algorithm with adaptive mesh refinement based on the damage length scale variable $\xi$}
\begin{algorithmic}[1]
\State \textbf{Input:} Initial mesh, model parameters, load increment \(\Delta u\), tolerance \(tol\)
\State \textbf{Input:} AMR parameters \(\xi_{\text{refine}}\), \(h_{\min}\); function \(g(\xi)\); flag \(\texttt{refinement\_enabled}\)
\State \textbf{Output:} Converged fields \(u\), \(v\)
\State Initialize fields \(u_0, v_0\) from initial conditions.
\State Set load step counter \(n = 0\).

\While{load step \(n < N_{\text{max}}\)} \Comment{Main loop over load increments}
    \State \textit{// Initialize for staggered iterations at current load step}
    \State Set staggered iteration counter \(k = 0\).
    \State Set initial guesses: \(u_{n,k} \gets u_{n-1}\), \(v_{n,k} \gets v_{n-1}\).
    \Repeat
        \State \textit{// Store previous iteration's solution for error calculation}
        \State \(u_{\text{prev}} \gets u_{n,k}\), \(v_{\text{prev}} \gets v_{n,k}\)
        \State \(k \gets k + 1\)

        \State \textbf{Solve for displacement \(u_{n,k}\)} using \(v_{\text{prev}}\).
        
        \State \textbf{Solve for phase-field \(v_{n,k}\)} using the newly computed \(u_{n,k}\).

        \State \textit{// Check for convergence of the staggered scheme}
        \State Compute $L_2$ norm of displacement error: \( \text{err}_u = \frac{\| u_{n,k} - u_{\text{prev}} \|_2}{\| u_{n,k} \|_2} \)
        \State Compute $L_2$ norm of phase-field error: \( \text{err}_v = \frac{\| v_{n,k} - v_{\text{prev}} \|_2}{\| v_{n,k} \|_2} \)
        
    \Until{\( (\text{err}_u < tol) \) \textbf{and} \( (\text{err}_v < tol) \)} \Comment{Staggered loop converges}

    \State \textit{// Solution at step n has converged: \(u_n \gets u_{n,k}\), \(v_n \gets v_{n,k}\)}
    
    \State \textit{// Adaptive Mesh Refinement based on converged phase-field solution}
    \If{\texttt{refinement\_enabled}}
        \State Set \texttt{refinement\_occurred} = \textbf{False}
        \For{each element \(e\) in the mesh}
            \State Let \(\xi_e\) be the value of the phase-field \(v_n\) in element \(e\).
            \State Let \(h\) be the current size of element \(e\).
            \If{$\xi_e < \xi_{\text{refine}}$ \textbf{and}  $h > h_{\min}$}
                \State Flag element \(e\) for refinement.
                \State Set \texttt{refinement\_occurred} = \textbf{True}
            \EndIf
        \EndFor
        \If{\texttt{refinement\_occurred}}
            \State Refine all flagged elements and project fields \(u_n, v_n\) onto the new mesh.
        \EndIf
    \EndIf

    \State \textit{// Prepare for the next load step}
    \State Update boundary conditions with the next load increment, e.g., \(u^{BC}_{n+1} \gets u^{BC}_{n} + \Delta u\).
    \State \(n \gets n + 1\).
\EndWhile
\end{algorithmic}
\end{algorithm}

\subsection{Problem Configuration }

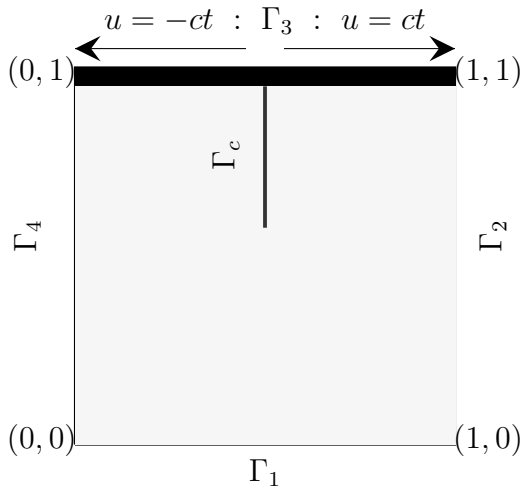
\begin{figure}[H]
\centering
\begin{tikzpicture}[scale=1.25]
    \tikzset{myptr/.style={decoration={markings,mark=at position 1 with %
    {\arrow[scale=3,>=stealth]{>}}},postaction={decorate}}}
    \filldraw[draw=black, thick] (0,0) -- (4,0) -- (4,4) -- (0,4) -- (0,0);
    \shade[inner color=gray!7, outer color=gray!7] (0,0) -- (4,0) -- (4,3.8) -- (0,3.8) -- (0,0);
    \node at (-0.35,0.05)   {$(0,0)$};
    \node at (4.35,0.05)   {$(1,0)$};
    \node at (-0.35,3.95)   {$(0,1)$};
    \node at (4.35, 3.95)   {$(1,1)$}; %~ \mathbf{n}\cdot \nabla u=0$
    \node at (-0.5, 2.5)[anchor=east, rotate=90]{$\Gamma_4$};
  %  \node at (-0.5, 3)[anchor=east, rotate=90]{$\&~ \boldsymbol{n} \cdot \nabla v=0$};
    \node at (4.4, 2.5)[anchor=east, rotate=90]{$\Gamma_2$};
   % \node at (4.8, 3)[anchor=east, rotate=90]{$\&~ \boldsymbol{n} \cdot \nabla v=0$};
    \node at (2, -0.3) {$\Gamma_{1}$};
    \node at (2, 4.5) {$ u=-ct ~: ~\Gamma_{3} ~:~ u=ct$};
   % \node at (2, 5) {$ ~ \boldsymbol{n} \cdot \nabla v=0$};
    \draw [myptr](1.8, 4.2)--(0.0, 4.2);
    \draw [myptr](2.2, 4.2)--(4.0, 4.2);
    %%%%%%%%%%%%%%%%
    \draw [line width=0.5mm, black!80]  (2,3.8) -- (2,2.3);
    \node at (1.59,3.35)[anchor=east, rotate=90]{$\Gamma_c$};
   % \node at (2.25,3.8)[anchor=east, rotate=90]{$\mathbf{n}\cdot \nabla u=0$};
   
\end{tikzpicture}
\caption{A computational domain showing an edge crack with anti-plane shear loading.}\label{fig:crackdomain}
\end{figure}

The numerical simulations are conducted on a square domain containing a pre-existing edge crack located at the center of the top boundary, as illustrated in Figure \ref{fig:crackdomain}. The material is subjected to an anti-plane shear loading, which induces a Mode III (tearing) fracture condition. This is achieved by applying time-dependent, opposing displacements to the top boundary, $\Gamma_3$, on either side of the crack:
\begin{align*}
u =
\begin{cases}
-c\,t  & \text{if } x \in (0, 0.5)  \quad \text{and} \quad y=1, \\
~~c\,t  & \text{if } x \in (0.5, 1) \quad \text{and} \quad y=1.
\end{cases}
\end{align*}
The boundary conditions are chosen to reflect a realistic physical scenario:
\begin{itemize}
    \item \textbf{Mechanical Boundary Conditions:} Outside of the loading region on $\Gamma_3$, all external surfaces, as well as the faces of the crack itself ($\Gamma_c$), are defined as traction-free ($\mathbf{n}\cdot \nabla u=0$). This means no external forces are applied to these surfaces, allowing them to deform freely according to the internal stress distribution.

    \item \textbf{Phase-Field Boundary Conditions:} A homogeneous Neumann condition ($\mathbf{n}\cdot \nabla v=0$) is enforced on all boundaries for the phase-field variable. This condition is fundamental to ensuring that the crack propagation is governed by the material's internal energy and not influenced by artificial constraints at the boundaries. Physically, it means that the boundaries are neutral; they neither initiate nor impede fracture. This allows the crack to propagate naturally, with its path and speed determined solely by the evolving stress field within the domain.
\end{itemize}
%%%%%%%%%%%%%%%%%

In the remainder of this paper, we analyze the influence of the phase-field length scale parameter, $\xi$, by presenting results from two distinct modeling frameworks. The first approach employs a conventional, globally constant $\xi$, whose value is optimized for the problem to achieve uniform regularization of the crack interface. The second, more advanced approach treats $\xi$ as an inhomogeneous field that varies locally, allowing the model to adapt the diffuse crack's thickness in response to the evolving stress and damage states. This comparison will highlight the benefits of an adaptive length scale in terms of accuracy and computational efficiency.

\subsection{Optimal and globally constant $\xi$}
In this section, we detail the numerical procedure for determining an optimal, globally constant regularization parameter, $\xi$. 
The approach involves solving the regularized phase-field model utilizing a modified version of Algorithm~1. 
A key modification from the algorithm is the exclusion of local mesh refinement; instead, a globally uniform mesh is maintained throughout the computation. 
At each time step, an optimal value for $\xi$ is calculated for the entire domain based on the numerical solution via Equation~\ref{eq:xi_const_xi}.

For the numerical simulations, the material properties were set to a fracture toughness of $G_c = 2.7$ and a shear modulus of $\mu=80.8$~GPa. 
All other parameters appearing in the model are adopted from the calibrated values listed in Table~\ref{tab:param_calib}.

\begin{figure}[htb!]
    \centering 
    % --- First subfigure ---
    \begin{subfigure}{0.3\textwidth}
        \centering
        \includegraphics[width=\linewidth]{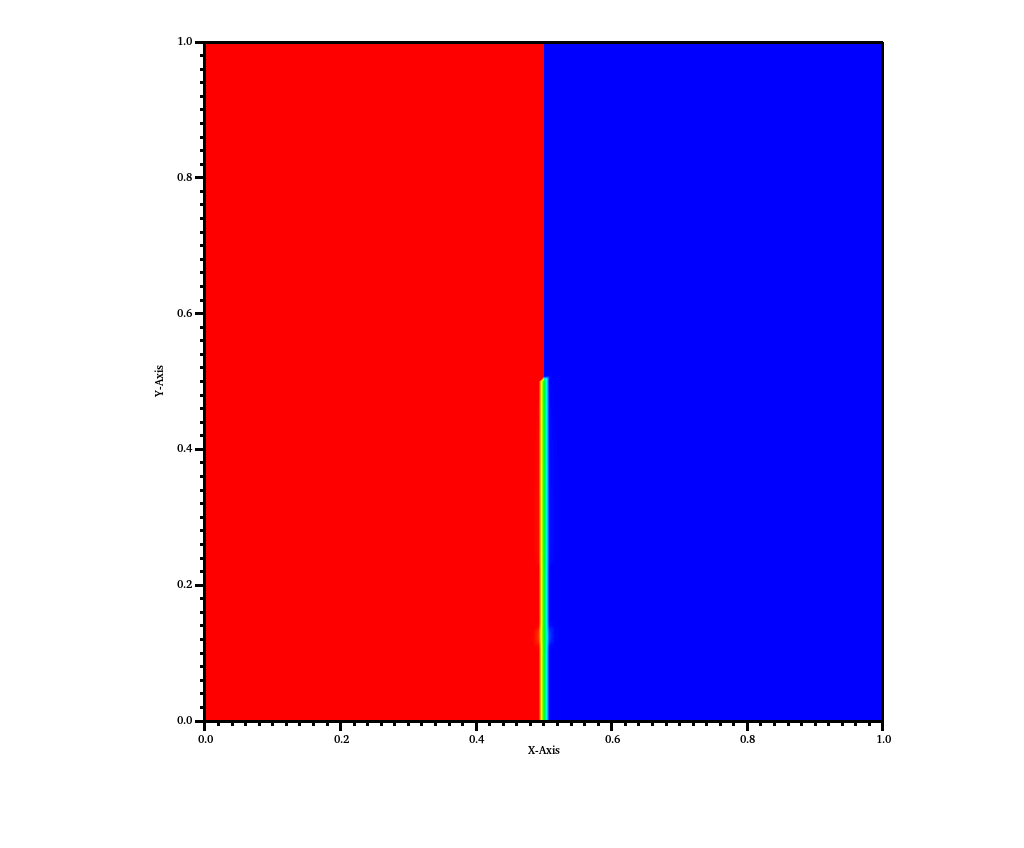}
        \caption{$128$ cells }
        \label{fig:sub1}
    \end{subfigure}
    \hfill 
    % --- Second subfigure ---
    \begin{subfigure}{0.3\textwidth}
        \centering
        \includegraphics[width=\linewidth]{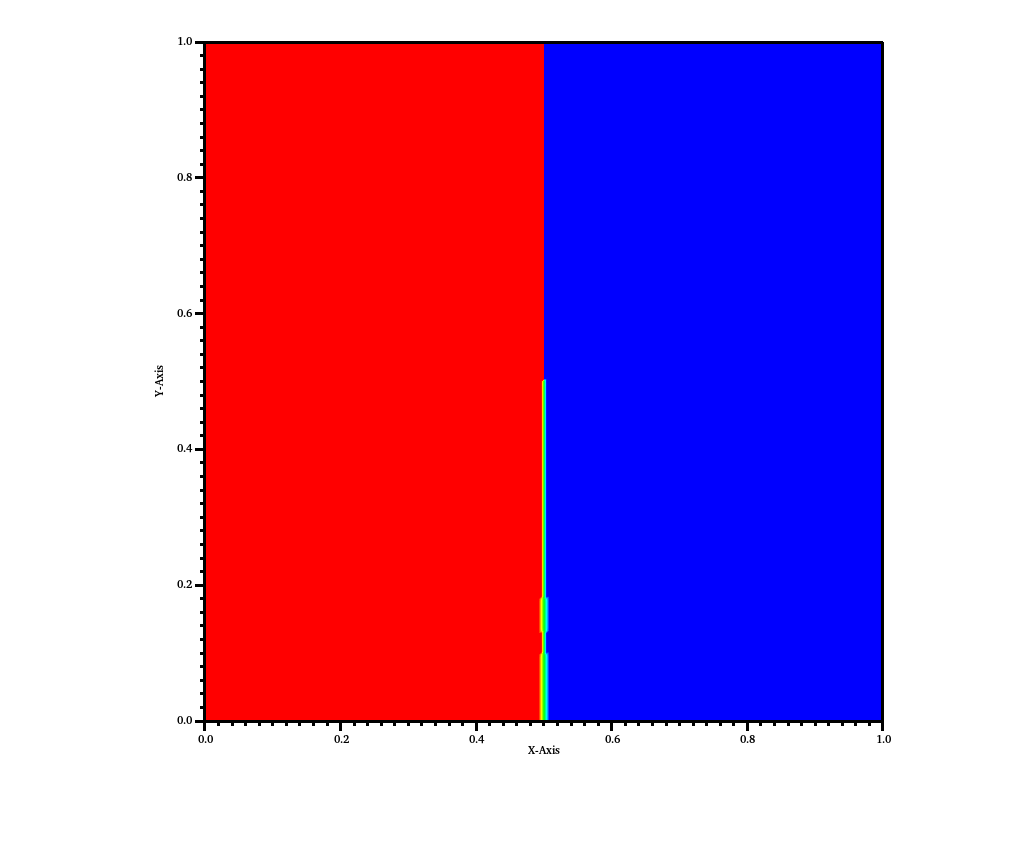}
        \caption{$256$ cells}
        \label{fig:sub2}
    \end{subfigure}
    \hfill 
    % --- Third subfigure ---
    \begin{subfigure}{0.3\textwidth}
        \centering
        \includegraphics[width=\linewidth]{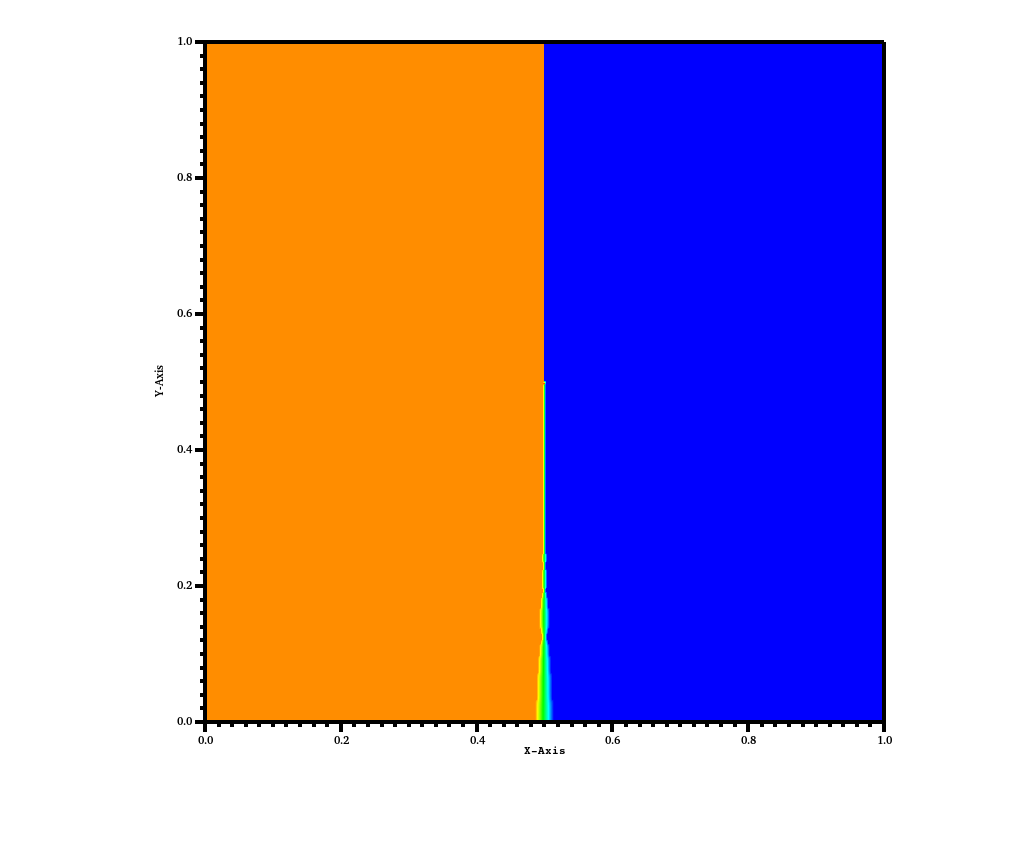}
        \caption{$512$ cells}
        \label{fig:sub3}
    \end{subfigure}
    %--------------------------------
        \vspace{1em} 
    % --- Second Row ---
        % --- First subfigure ---
    \begin{subfigure}{0.3\textwidth}
        \centering
        \includegraphics[width=\linewidth]{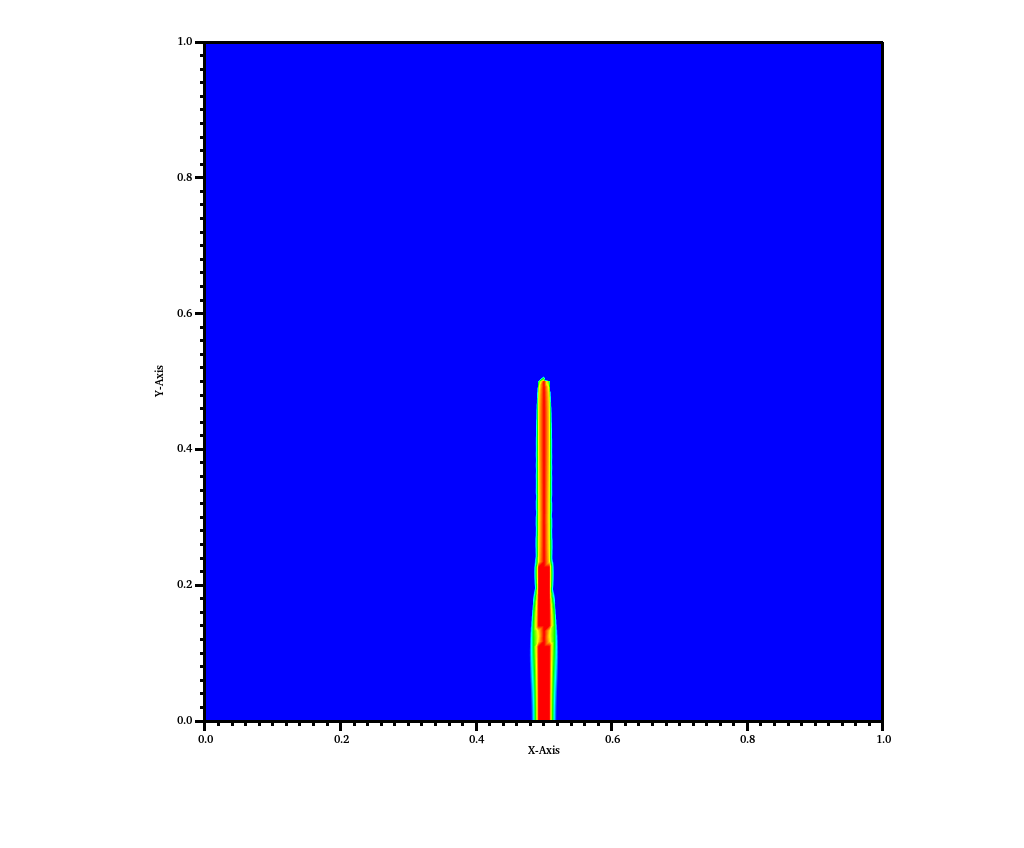}
        \caption{$128$ cells}
        \label{fig:sub1}
    \end{subfigure}
    \hfill 
    % --- Second subfigure ---
    \begin{subfigure}{0.3\textwidth}
        \centering
        \includegraphics[width=\linewidth]{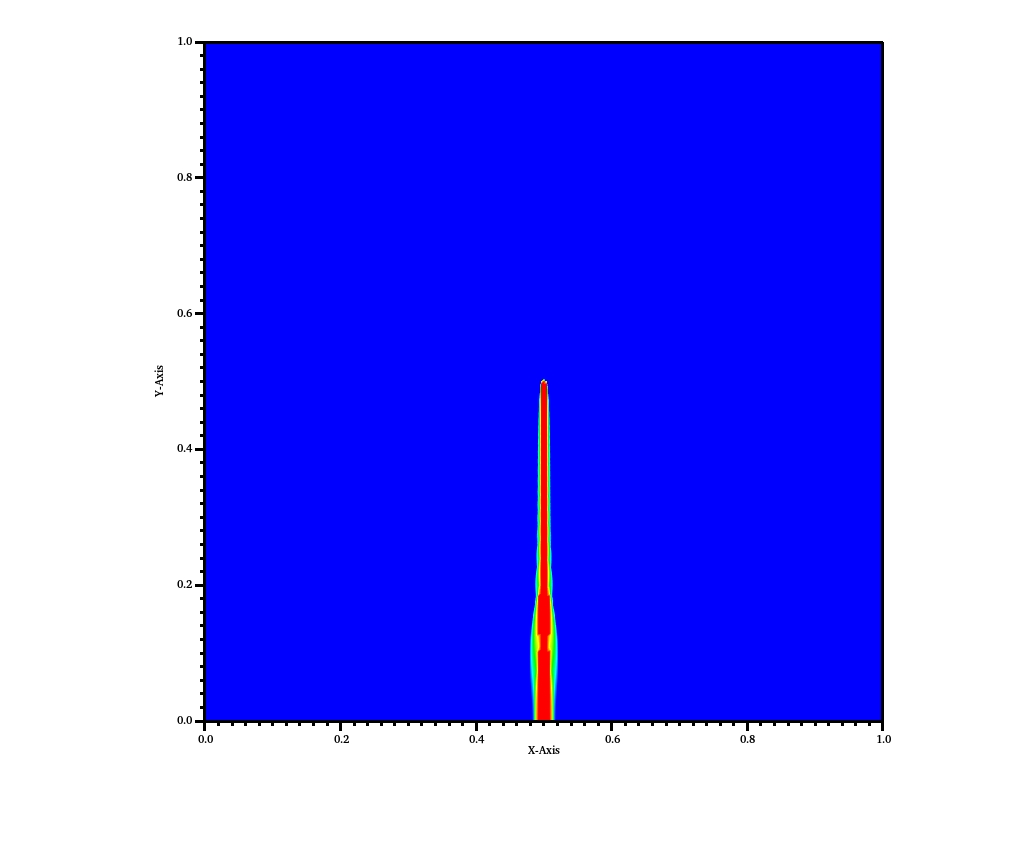}
        \caption{$256$ cells}
        \label{fig:sub2}
    \end{subfigure}
    \hfill 
    % --- Third subfigure ---
    \begin{subfigure}{0.3\textwidth}
        \centering
        \includegraphics[width=\linewidth]{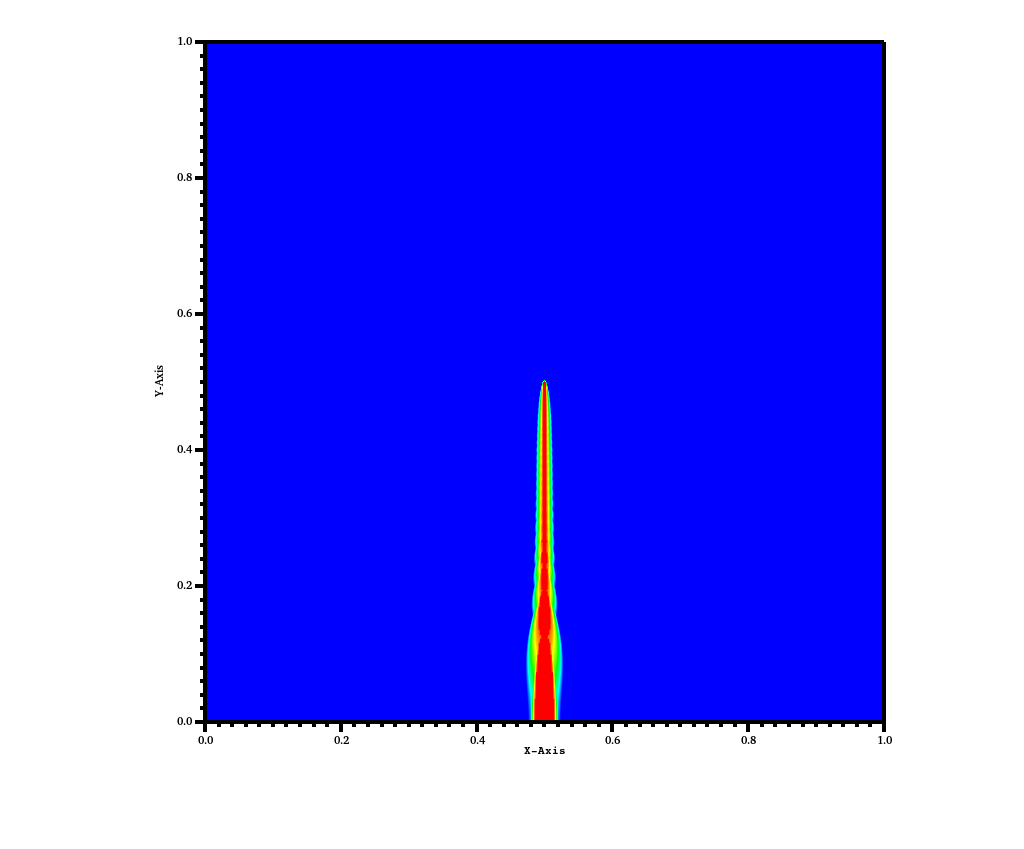}
        \caption{$512$ cells}
        \label{fig:sub3}
    \end{subfigure}
    \caption{Plof of $u$ and $v$ on the different meshes with globally constant $\xi$.}
    \label{u_v_globalConsXI}
\end{figure}

Figure~\ref{u_v_globalConsXI} presents the numerical results for the displacement field, $u$, and the phase-field variable, $v$, obtained using various mesh resolutions. 
The contour plots for the phase-field variable, $v$, illustrate the predicted crack topology. 
In these plots, the red regions, where $v \approx 0$, represent the fully fractured material, while the blue regions, where $v \approx 1$, correspond to the undamaged material state. 
Concurrently, the plots for $u$ depict the out-of-plane displacement field, which is characteristic of the applied anti-plane shear loading conditions.

\begin{table}[h!]
    \centering
    \caption{The optimal value of $\xi$ computed from Equation~\ref{eq:xi_const_xi}.}
    \label{tab:xi}
    \begin{tabular}{cccc}
        \toprule
        {n} & {h} & {${\xi = m \cdot h}$} & {optimal ${\xi}$} \\
        \midrule
        128 & 0.008 & 0.0390625 & 0.13687 \\
        256 & 0.004 & 0.03125   & 0.06927 \\
        512 & 0.002 & 0.0234375   & 0.03464 \\
        \bottomrule
    \end{tabular}
\end{table}

Table~\ref{tab:xi} summarizes the results of our analysis, comparing the optimal regularization parameter, $\xi$, computed via Equation~\ref{eq:xi_const_xi}, against the smaller, pre-selected value used at the start of the minimization process. 
For each mesh resolution tested, the optimally determined $\xi$ is notably larger. 
This outcome is a significant advantage, as a larger regularization parameter relaxes the system's stiffness, thereby accelerating computational performance.

To validate this approach, we also performed simulations using the smaller, conventional $\xi$ values. 
These tests yielded nearly identical crack profiles and displacement fields, confirming that the improved computational efficiency does not compromise the accuracy of the results. 
Therefore, the proposed method provides an attractive and systematic framework for selecting the influential $\xi$ parameter in the \textsf{AT1} model, resulting in substantial savings in computational time and resources without compromising predictive accuracy. 

\begin{figure}[htb!]
    \centering
    % --- First subfigure ---
    \begin{subfigure}{0.32\textwidth}
        \centering
        \includegraphics[width=\linewidth]{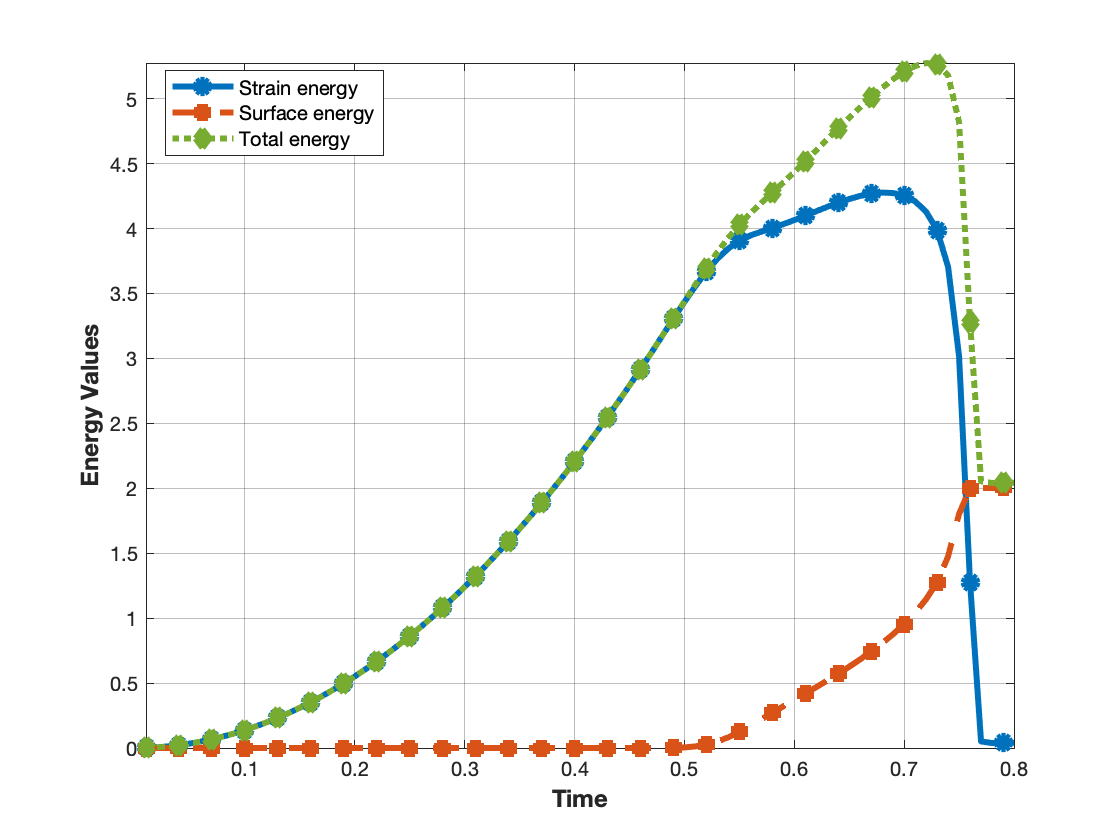}
        \caption{$128$ cells}
        \label{fig:sub1}
    \end{subfigure}\hfill
    % --- Second subfigure ---
    \begin{subfigure}{0.32\textwidth}
        \centering
        \includegraphics[width=\linewidth]{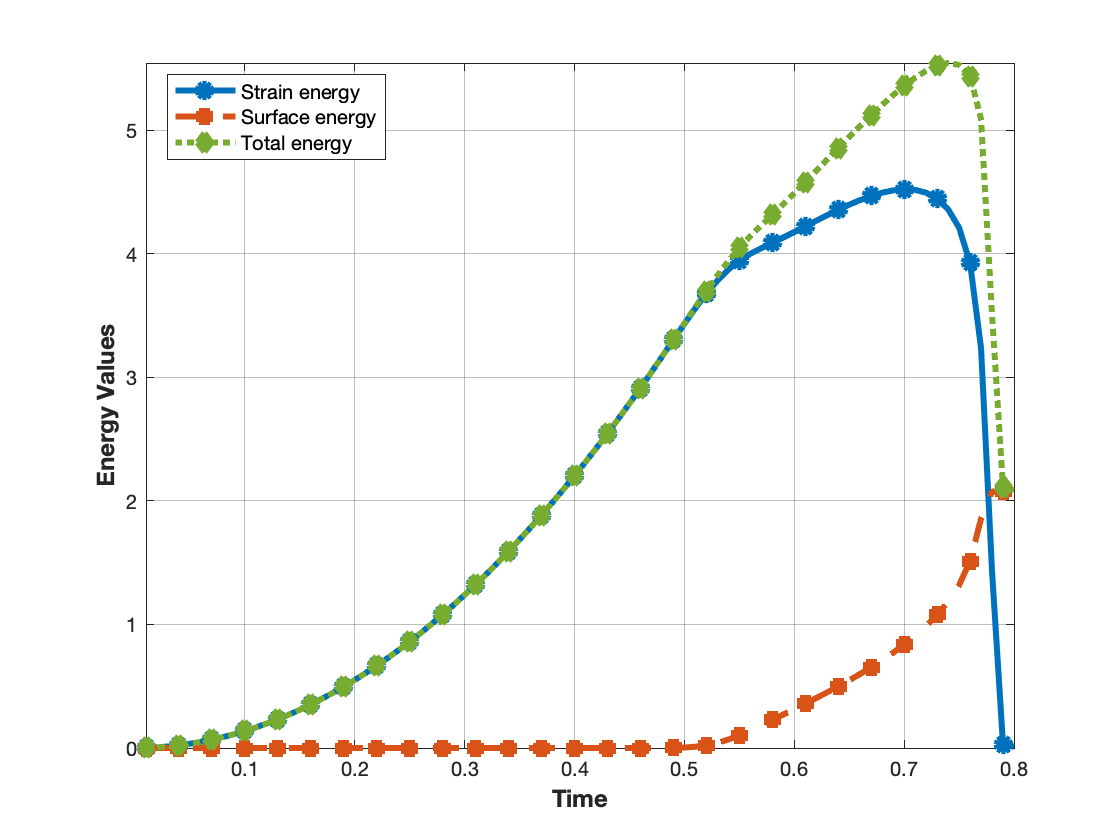}
        \caption{$256$ cells}
        \label{fig:sub2}
    \end{subfigure}\hfill
    % --- Third subfigure ---
    \begin{subfigure}{0.32\textwidth}
        \centering
        \includegraphics[width=\linewidth]{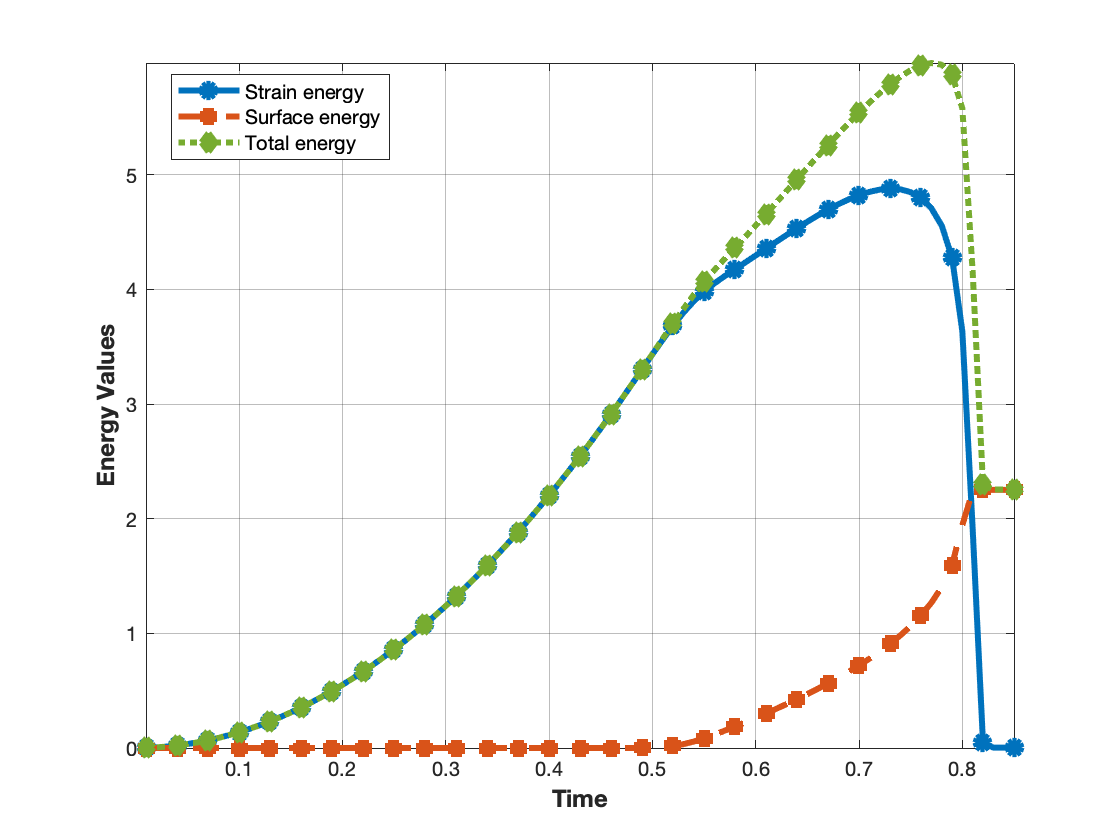}
        \caption{$512$ cells}
        \label{fig:sub3}
    \end{subfigure}
  \caption{Evolution of energies (strain, surface, and total) from an \textsf{AT1} phase-field fracture simulation for different computational mesh resolutions. The plots show the conversion of stored strain energy into surface energy as the crack propagates, with results converging as the mesh is refined.}
    \label{glob_mesh_energies}
\end{figure}

The Figure~\ref{glob_mesh_energies}  illustrate the evolution of strain, surface, and total energies over time from a brittle fracture simulation using the Ambrosio-Tortorelli ( \textsf{AT1}) phase-field model. Each plot corresponds to a different computational mesh resolution ( $128$, $256$, and $512$ cells), demonstrating the model's behavior and numerical convergence. Initially, as the system is loaded, the strain energy (blue line) increases, representing the storage of elastic energy in the material. The surface energy (orange line), which approximates the energy required to create a crack, remains near zero. At a critical point around time$\approx 0.7$, the stored strain energy is abruptly released, causing a catastrophic failure. This release of strain energy drives the creation of a crack, marked by the simultaneous and rapid increase in surface energy. The total energy (green line), representing the sum of the two, illustrates the overall energy balance, which peaks just before the final fracture event. Comparing the three plots reveals excellent convergence with mesh refinement; the overall behavior is consistent, with only slight increases in the peak strain energy and a minor delay in the fracture time as the mesh becomes finer. This behavior is characteristic of how phase-field models capture the sharp transition from energy storage to crack propagation.

\subsection{Spatially varying (or an inhomogeneous) $\xi$}
In this section, we solve the minimization problem for displacement $u$, phase field $v$, and a globally varying damage length scale parameter $\xi$. The Algorithm-1 is implemented using a coarse mesh with $64$-cells, and then the AMR strategy using dynamically varying $\xi$ is implemented with an upper bound of $4$-refinement for a cell (flagged for refinement). Hence, the mesh size has a clear lower and upper bound, i.e. $1/2^6 \leq h \leq 1/2^{10}$. This yields bounds for $\xi$ also, i.e., $0.15 \leq \xi \leq 0.011$. The AMR strategy also includes coarsening the cells that were flagged with a mesh size smaller than the smallest prescribed mesh size; hence, the AMR strategy will not become trapped in an infinite refinement loop and will terminate after a finite number of refinements. 

For the numerical simulation, the material properties were defined by a fracture toughness of $G_c = 2.7$ and a shear modulus of $\mu=80.8$ GPa. The remaining model parameters were set to $\zeta = 31640.6$, $\alpha = 3.125$, and $\eta = 1.0 \times 10^{-10}$. The simulation proceeded with a fixed time step of $0.01$ and was run until the phase field, $v$, representing the crack, traversed the entire domain. A key aspect of this study involves the length scale parameter, $\xi$. Conventionally, $\xi$ is chosen to be proportional to the mesh size, $h$. For instance, with an initial coarse mesh size of $h=0.015625$, the corresponding parameter would be $\xi = 10h = 0.15625$. For the finest mesh used ($h=0.00097656$), this rule would imply a much smaller $\xi = 0.0097656$. However, our simulation dynamically determined the optimal value of $\xi$. At the first time step, the optimal value was found to be $\xi = 0.03464$. As the crack began to propagate (i.e., as $v$ decreased from its initial value of 1), the optimal $\xi$ settled into a range of $0.02419$ to $0.03644$. It is noteworthy that these computationally-derived optimal values are significantly larger than the value of $0.0097656$ that would typically be used for the fine mesh in other studies, highlighting a key departure from conventional practice \cite{yoon2021quasi}.
\begin{figure}[H]
    \centering
    \begin{subfigure}[b]{0.48\textwidth}
        \centering
        \includegraphics[width=\textwidth]{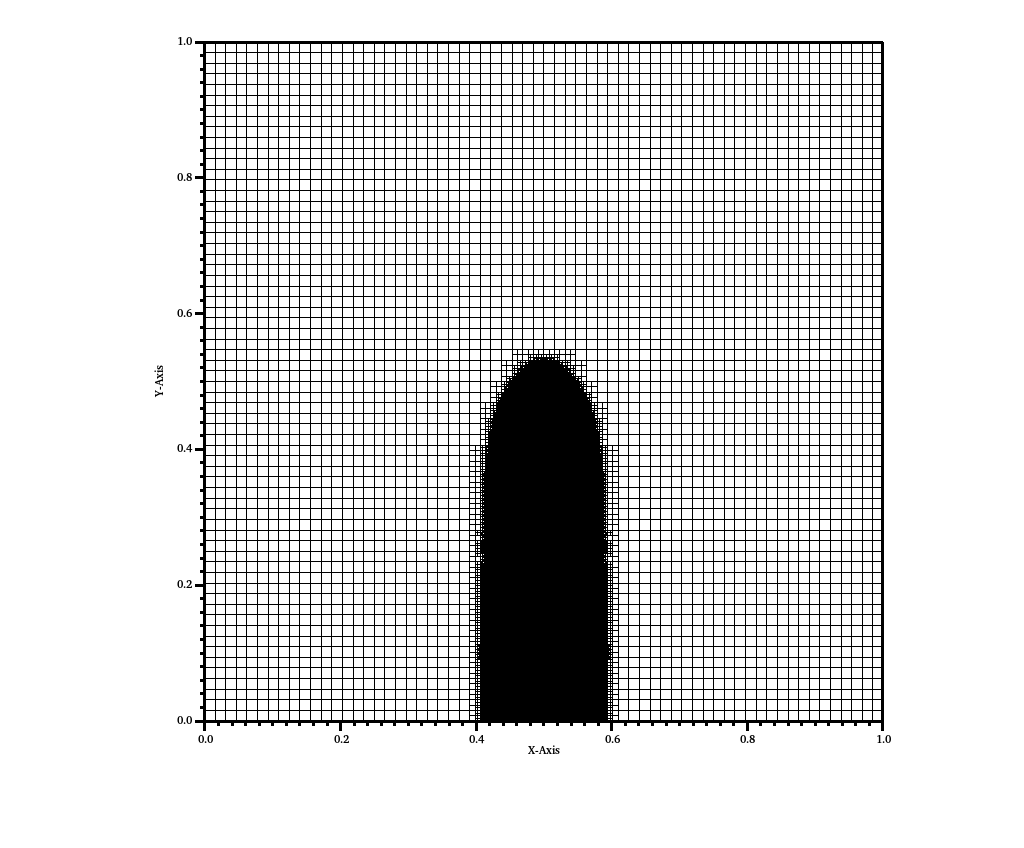}
        \caption{Mesh}
        \label{fig:mesh}
    \end{subfigure}
    \hfill 
    \begin{subfigure}[b]{0.48\textwidth}
        \centering
        \includegraphics[width=\textwidth]{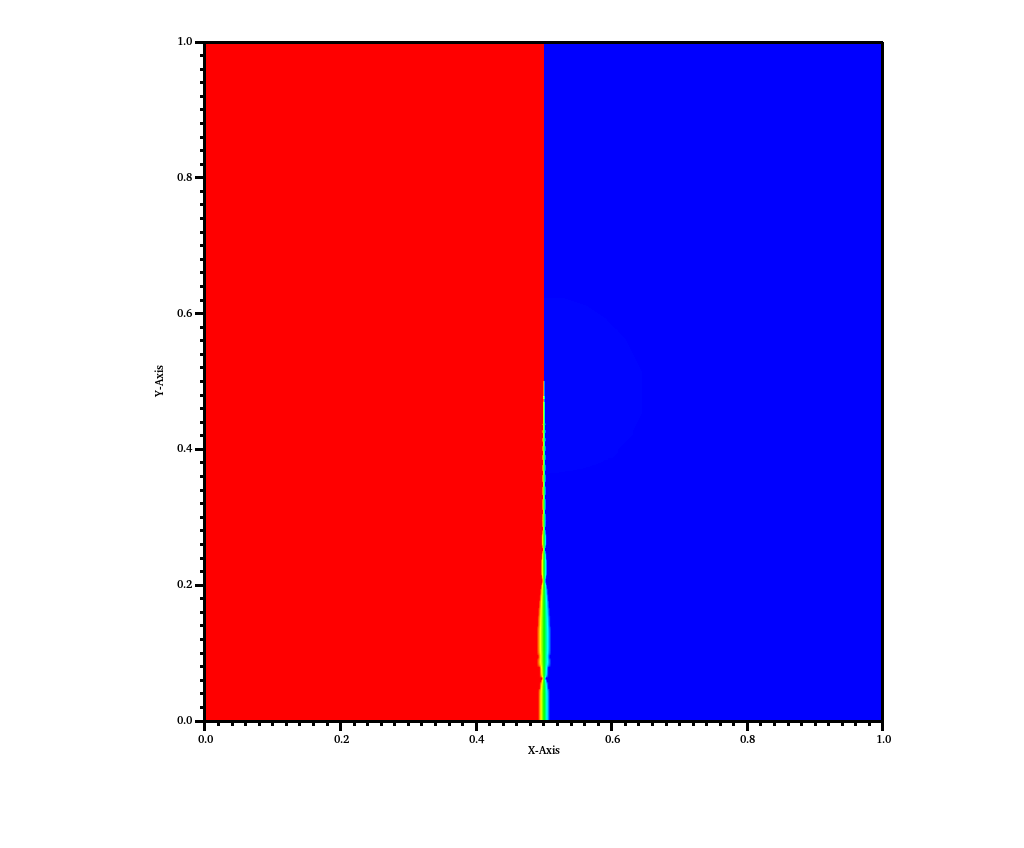}
        \caption{Displacement $u$}
        \label{fig:displacement}
    \end{subfigure}
    %-------------------
 \vspace{1em} % Adds some vertical space between rows
    % --- Second Row ---
        \begin{subfigure}[b]{0.48\textwidth}
        \centering
        \includegraphics[width=\textwidth]{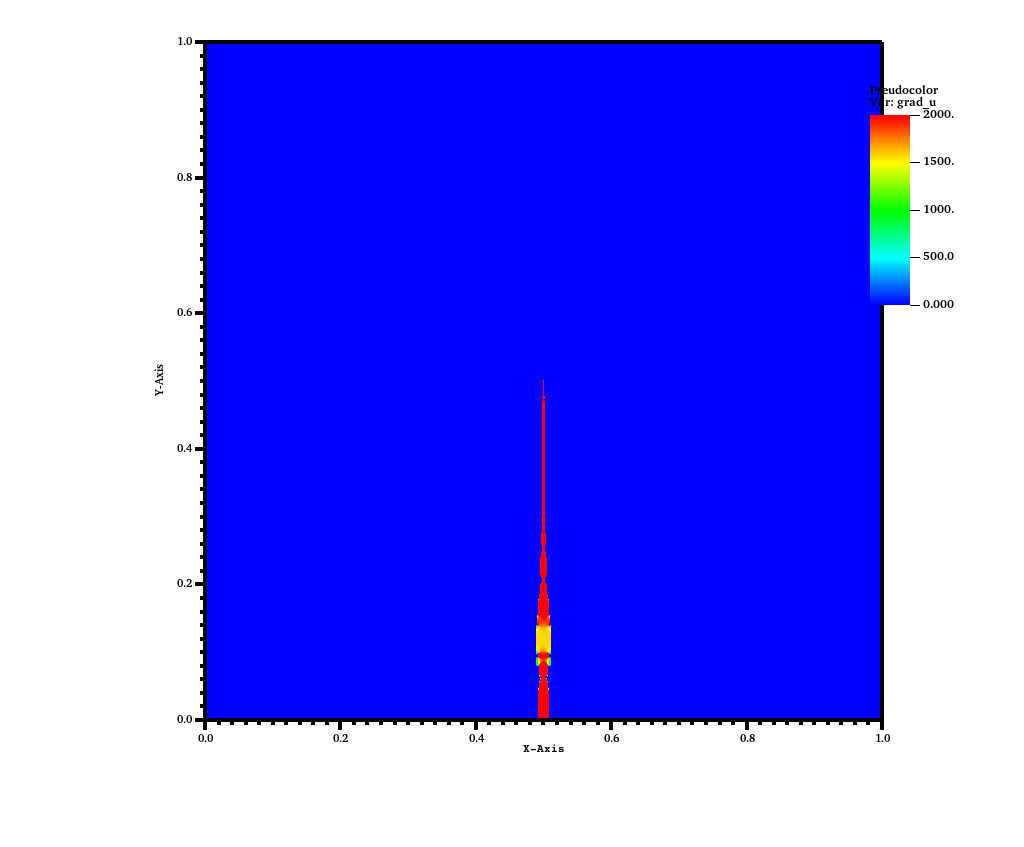}
        \caption{$\| \nabla u \|$}
        \label{fig:phase_field}
    \end{subfigure}
        \hfill 
    \begin{subfigure}[b]{0.48\textwidth}
        \centering
        \includegraphics[width=\textwidth]{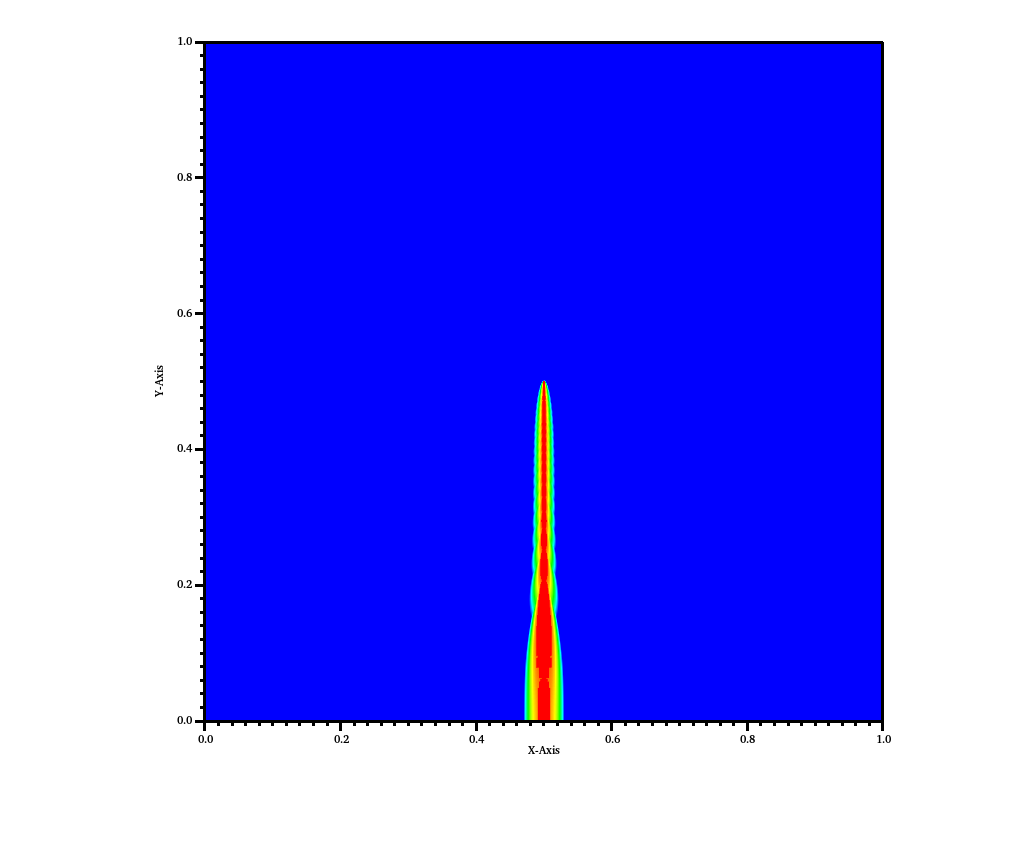}
        \caption{Phase field $v$}
        \label{fig:phase_field}
    \end{subfigure}
\caption{Top row: Computational mesh and the anti-plane displacement field $u$. Bottom row: The gradient magnitude, $\| \nabla u\|$, and the phase field, $v$. These simulation results correspond to a locally varying $\xi$.}
    \label{fig:mesh_u_v_LV_XI}
\end{figure}

Figure~\ref{fig:mesh_u_v_LV_XI} presents the key results from the simulation featuring a locally varying length-scale parameter, $\xi$. The top-left panel displays the computational mesh, which has been locally refined through an adaptive strategy. This refinement is driven by the value of $\xi$, computed using Equation~\eqref{eq:xi_field_xi}, ensuring high resolution in the regions where the phase field $v < 1$. This targeted adaptation yields a solution far superior to what would be achievable on a uniformly coarse mesh. The top-right panel shows the anti-plane displacement field, $u$, where the contours clearly illustrate the shear deformation characteristic of this mode of fracture. The bottom-left panel depicts the magnitude of the displacement gradient, $\| \nabla u \|$, which represents the elastic stress field. As expected, the stress is highly concentrated, with its maximum values occurring precisely along the path of the crack. Finally, the bottom-right panel shows the phase field, $v$, at the final time step. This field visualizes the fully developed crack path, where the dark red region ($v=0$) indicates the completely fractured material and the blue region ($v=1$) represents the undamaged material. 

\begin{figure}[H]
    \centering
    % --- First Row ---
    \begin{subfigure}[b]{0.48\textwidth}
        \centering
        \includegraphics[width=\textwidth]{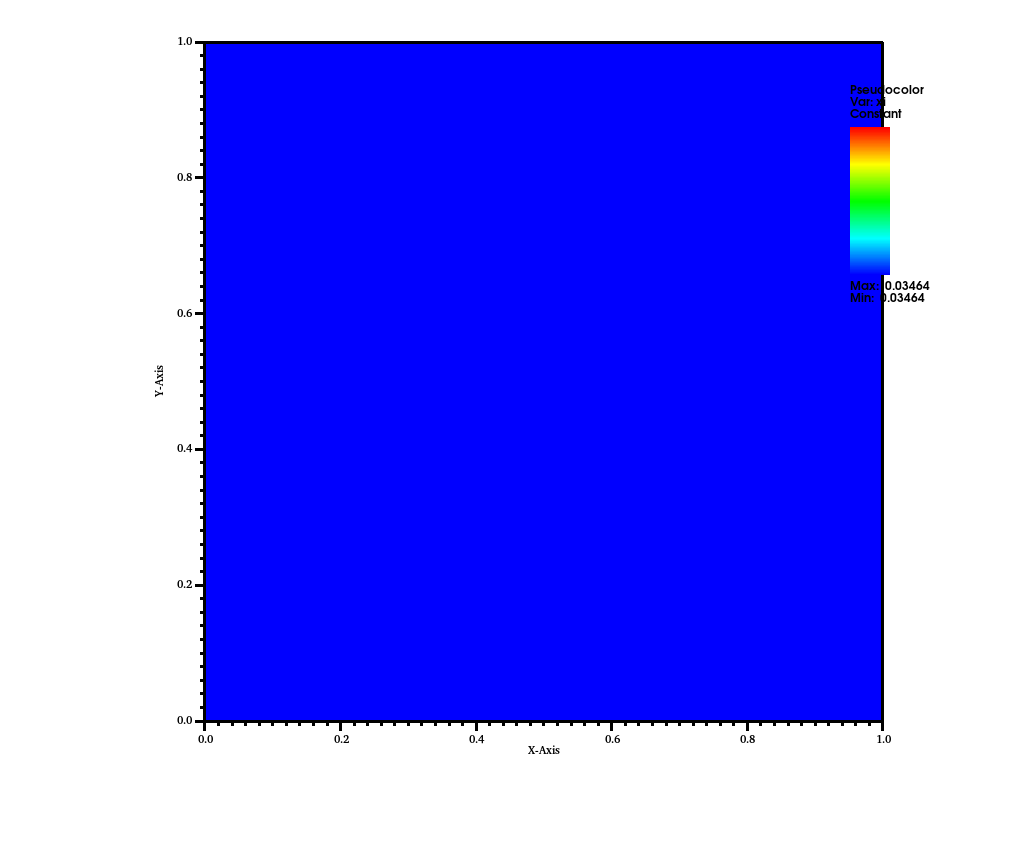}
        \caption{$n=1$}
        \label{fig:xi_n1}
    \end{subfigure}
    \hfill 
    \begin{subfigure}[b]{0.48\textwidth}
        \centering
        \includegraphics[width=\textwidth]{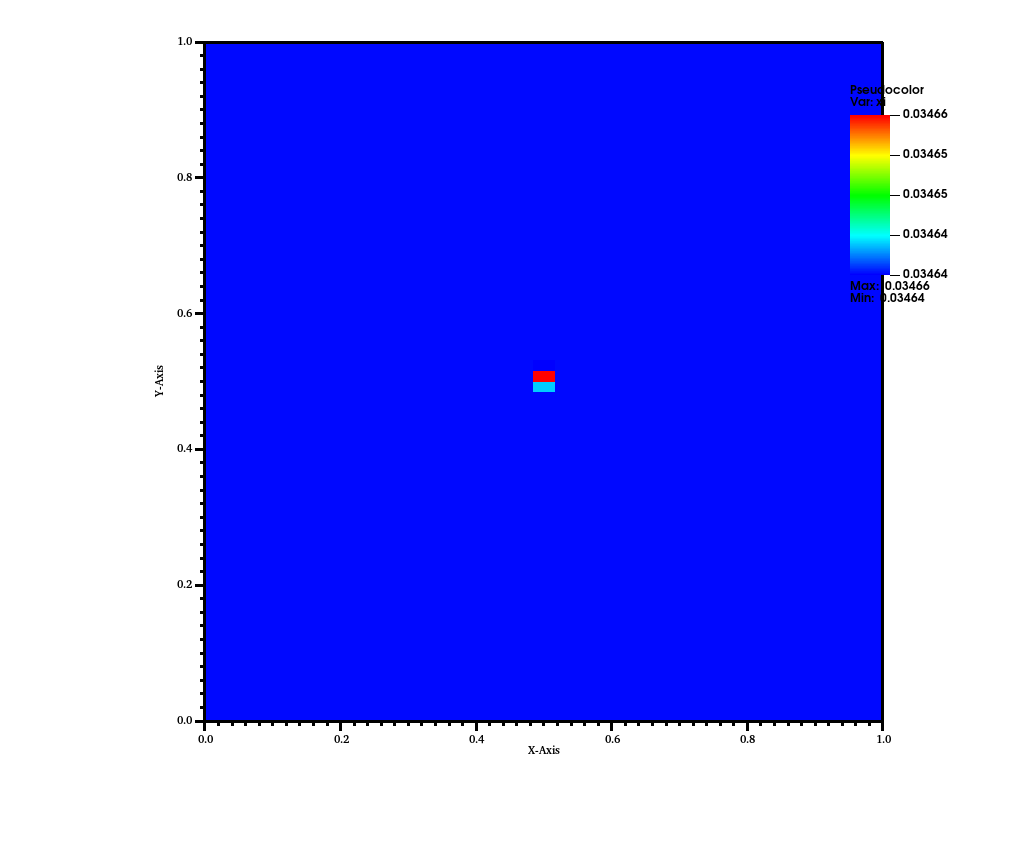}
        \caption{$n=50$}
        \label{fig:xi_n50}
    \end{subfigure}
    
    \vspace{1em} % Adds some vertical space between rows
    % --- Second Row ---
    \begin{subfigure}[b]{0.48\textwidth}
        \centering
        \includegraphics[width=\textwidth]{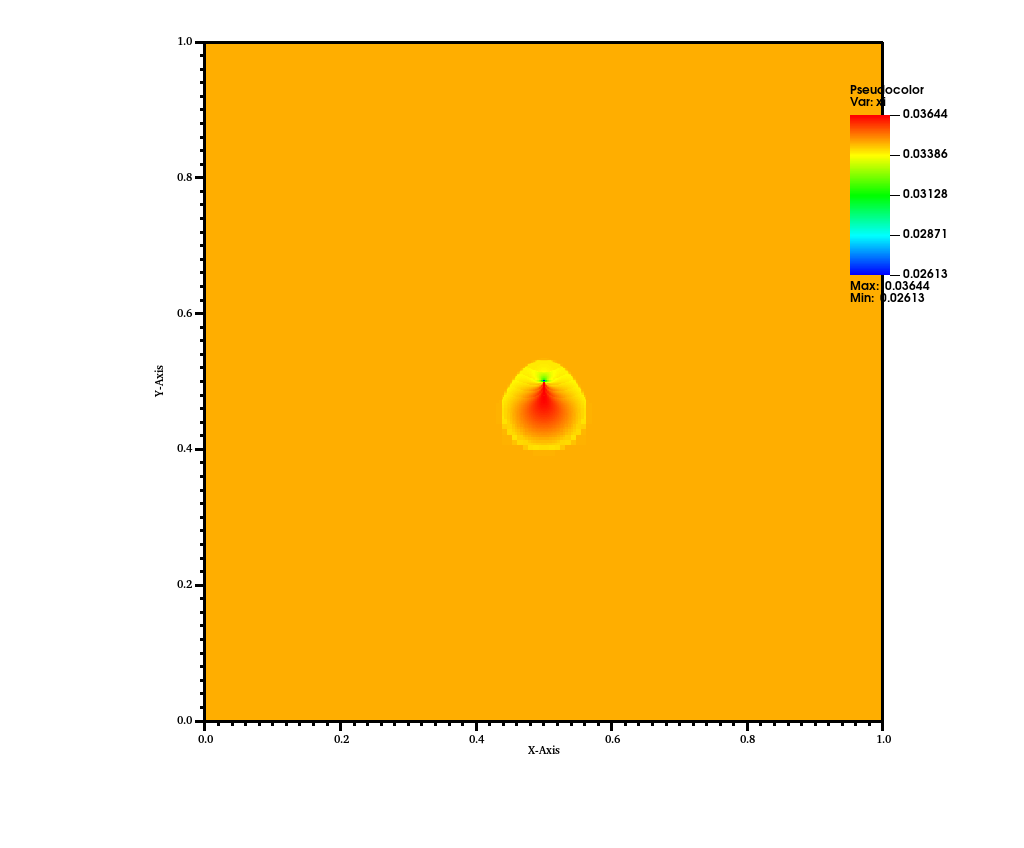}
        \caption{$n=60$}
        \label{fig:xi_n60}
    \end{subfigure}
    \hfill
    \begin{subfigure}[b]{0.48\textwidth}
        \centering
        \includegraphics[width=\textwidth]{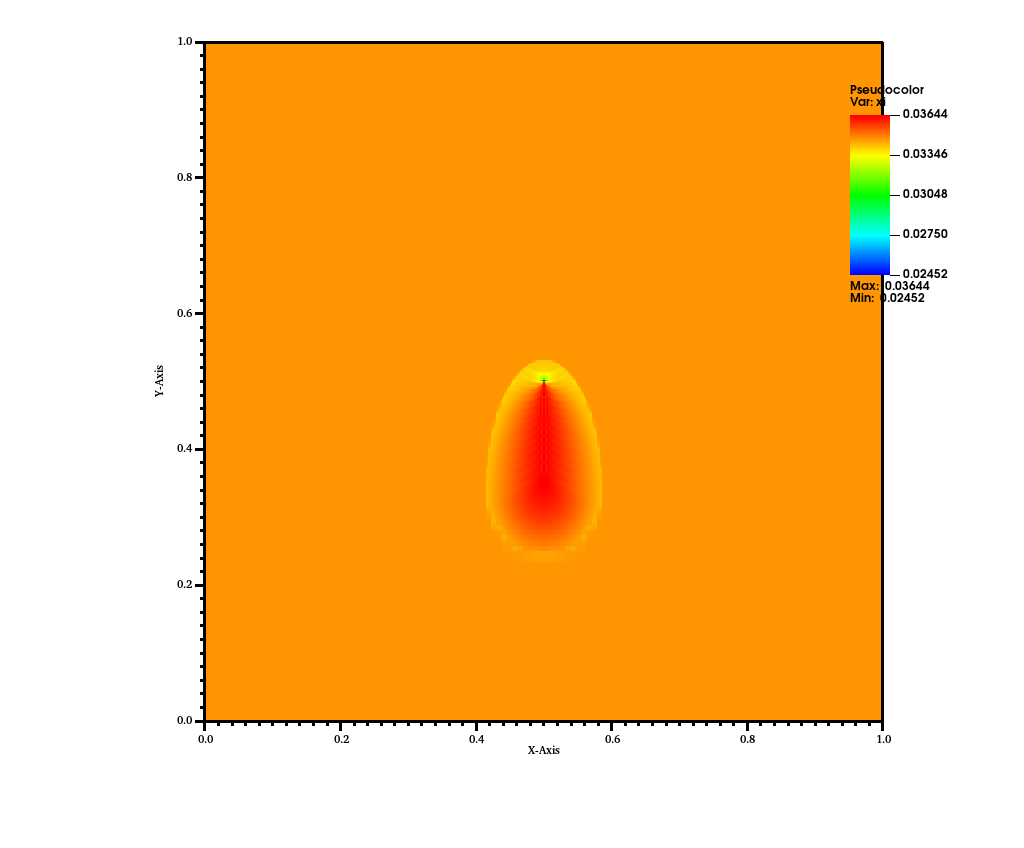}
        \caption{$n=70$}
        \label{fig:xi_n70}
    \end{subfigure}
    
    \vspace{1em} % Adds some vertical space between rows
    % --- Third Row ---
    \begin{subfigure}[b]{0.48\textwidth}
        \centering
        \includegraphics[width=\textwidth]{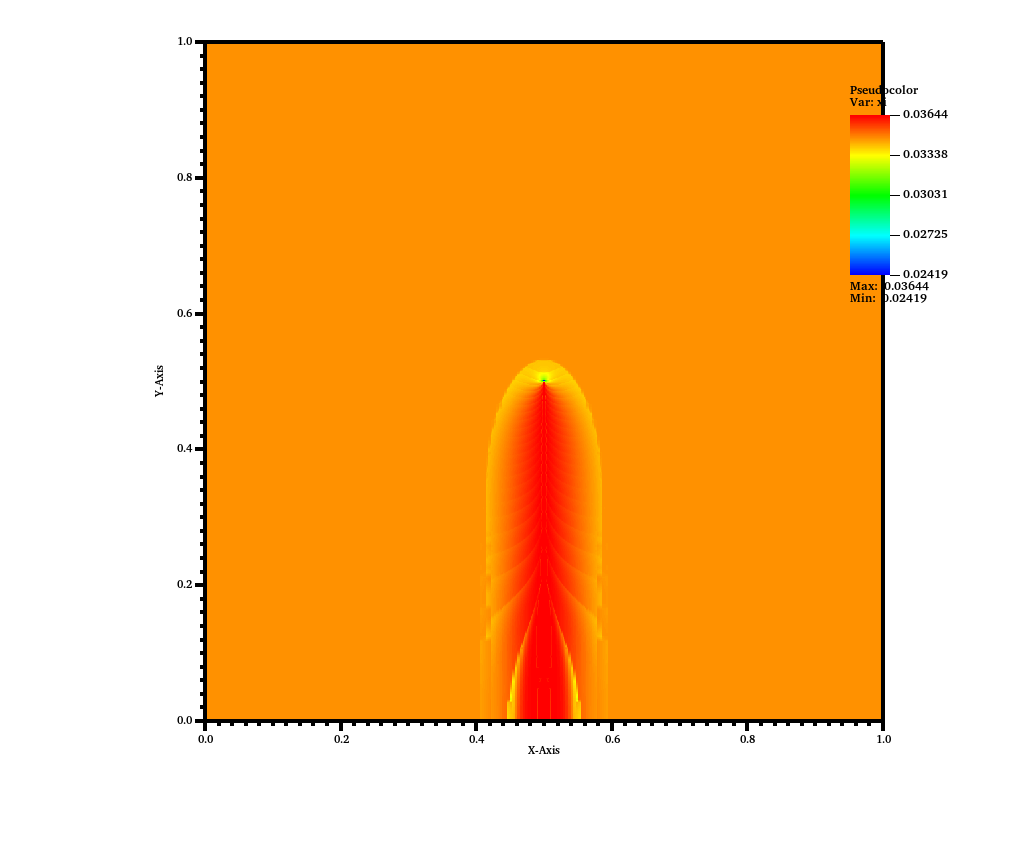}
        \caption{$n=80$}
        \label{fig:xi_n80}
    \end{subfigure}
    \hfill 
    \begin{subfigure}[b]{0.48\textwidth}
        \centering
        \includegraphics[width=\textwidth]{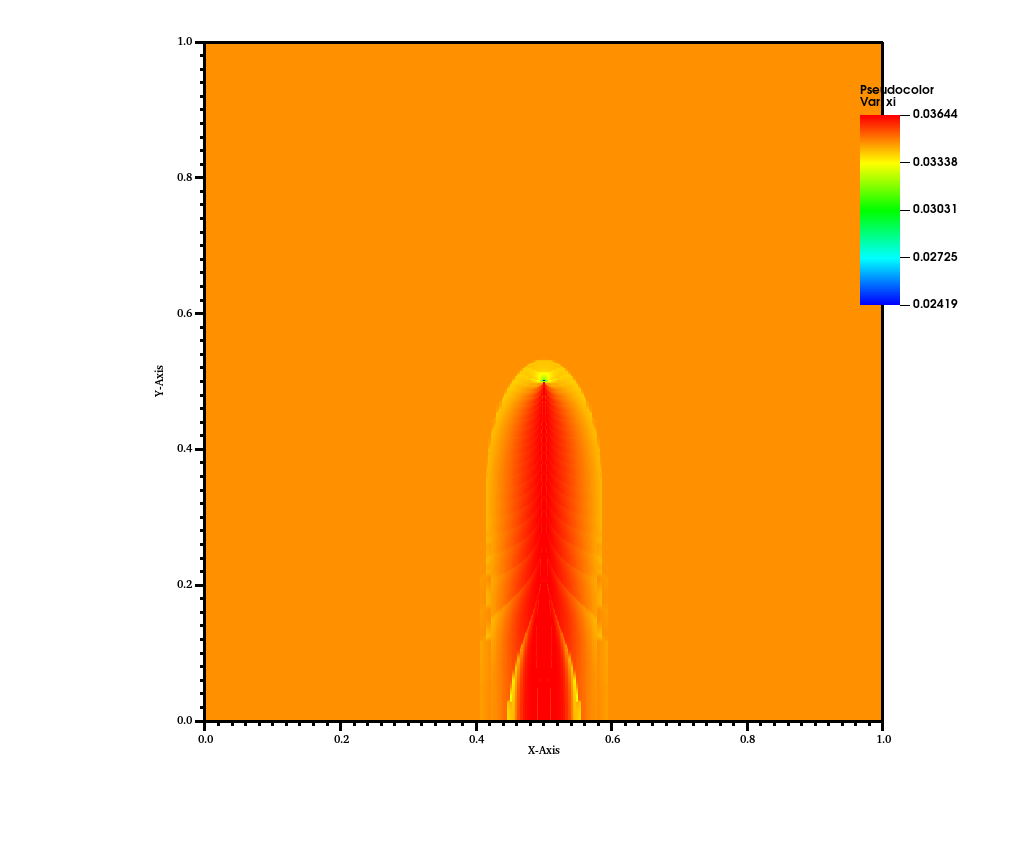}
        \caption{$n=85$}
        \label{fig:xi_n85}
    \end{subfigure}
    \caption{The evolution of $\xi$ at different time-steps ($n$).}
    \label{fig:xi_evolution}
\end{figure}

Figure~\ref{fig:xi_evolution} illustrates the dynamic evolution of the locally varying parameter $\xi$. 
Initially, for approximately the first 50 time steps, the value of $\xi$ remains constant at $0.03464$. 
This initial value is substantially larger than the threshold of $10h_f$, where $h_f$ represents the finest mesh size for a grid of $1024$ cells.  Subsequently, as the variable $v$ begins to decrease from its initial value of $1$, it prompts a corresponding decrease in $\xi$.  By the 85th time step, $\xi$ has reached its minimum value of $0.02419$ in this simulation. 
Notably, this minimum value is still an order of magnitude greater than the typical globally constant values used for $\xi$. 
This demonstrates a significant advantage of employing a locally varying $\xi$ within the \textsf{AT1} model, as it allows for adaptive and more physically representative parameterization.

\begin{figure}[H]
    \centering
    \includegraphics[width=0.6\textwidth]{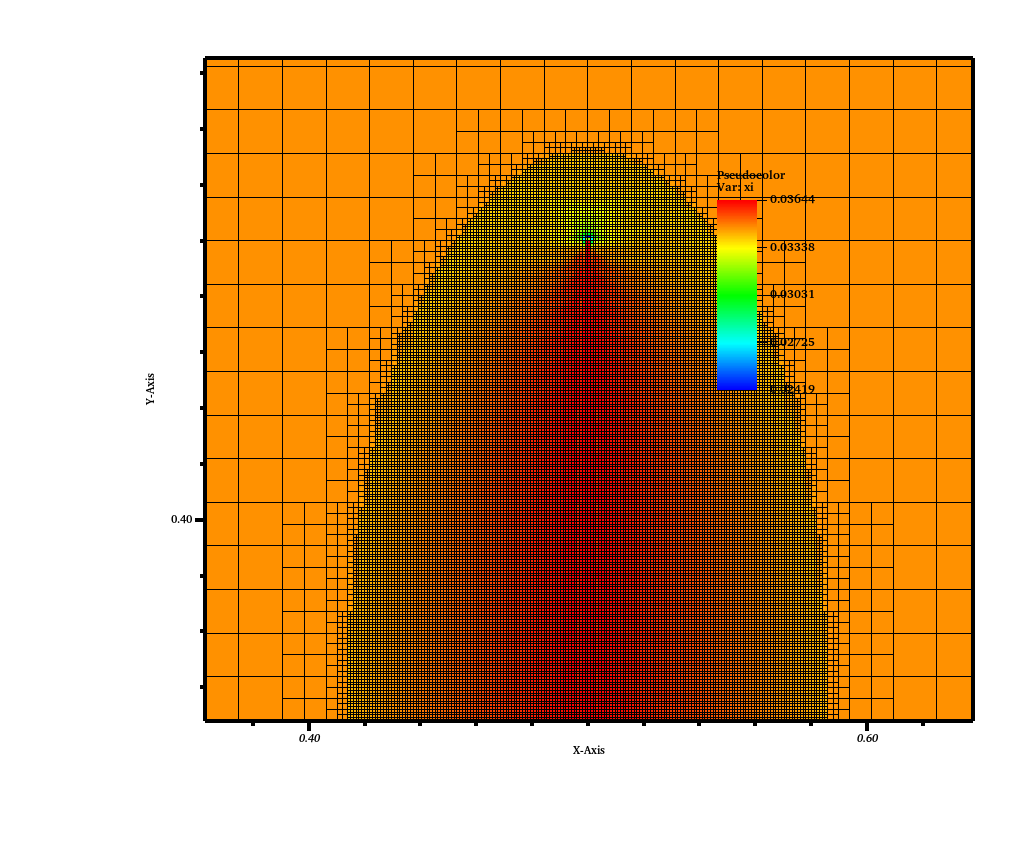}
\caption{The spatially varying $\xi$ near the crack-tip.}
    \label{xi_zoom}
\end{figure}
 
Figure~\ref{xi_zoom} provides a magnified view of the dynamically optimized length scale parameter, $\xi$, plotted along the trajectory of the fracture. The results clearly illustrate the adaptive nature of the parameter: it attains its minimum value precisely at the crack tip, where a smaller length scale is required to resolve the high-gradient damage zone. Conversely, moving away from the fractured region into the bulk material, the value of $\xi$ smoothly increases, approaching its maximum in the undamaged areas.

\begin{figure}[H]
    \centering
    \includegraphics[width=0.6\textwidth]{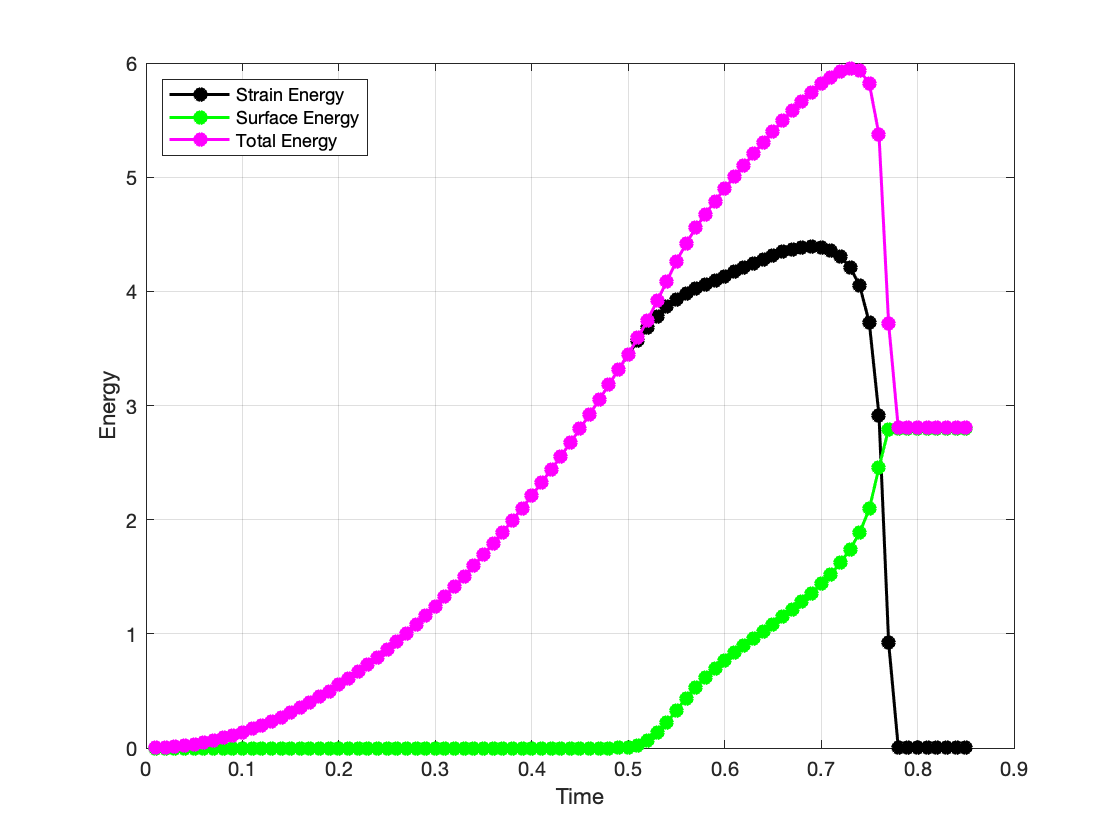}
\caption{A comparison of the energetic components (strain, surface, and total) for a model incorporating a locally varying length scale parameter, $\xi$.}
    \label{fig:energies_LA_XI}
\end{figure}

Figure~\ref{fig:energies_LA_XI} illustrates the classic energy dynamics of material fracture from the phase-field simulation, demonstrating the conversion of stored elastic energy into surface energy. Initially, during the elastic loading phase from $t \approx 0$ to $0.5$, the system stores energy purely through deformation. In this stage, the strain energy and total energy are identical and increase together, while the surface energy remains zero as no damage has yet occurred. The onset of damage begins at $t \approx 0.5$, marked by the rise in surface energy, after which the input energy is partitioned between stored strain energy and the energy dissipated to create new crack surfaces. This period of stable growth culminates in catastrophic failure as the strain energy peaks at $t \approx 0.75$, representing the material's critical load capacity. Subsequently, the stored elastic energy is rapidly converted into surface energy, causing the strain energy to plummet while driving swift crack propagation. This event leads to the final post-failure state for $t > 0.8$, in which the strain energy is fully dissipated, and the system's total energy becomes constant, equaling the final surface energy of the newly formed crack.

\begin{figure}[H]
    \centering
    \includegraphics[width=0.6\textwidth]{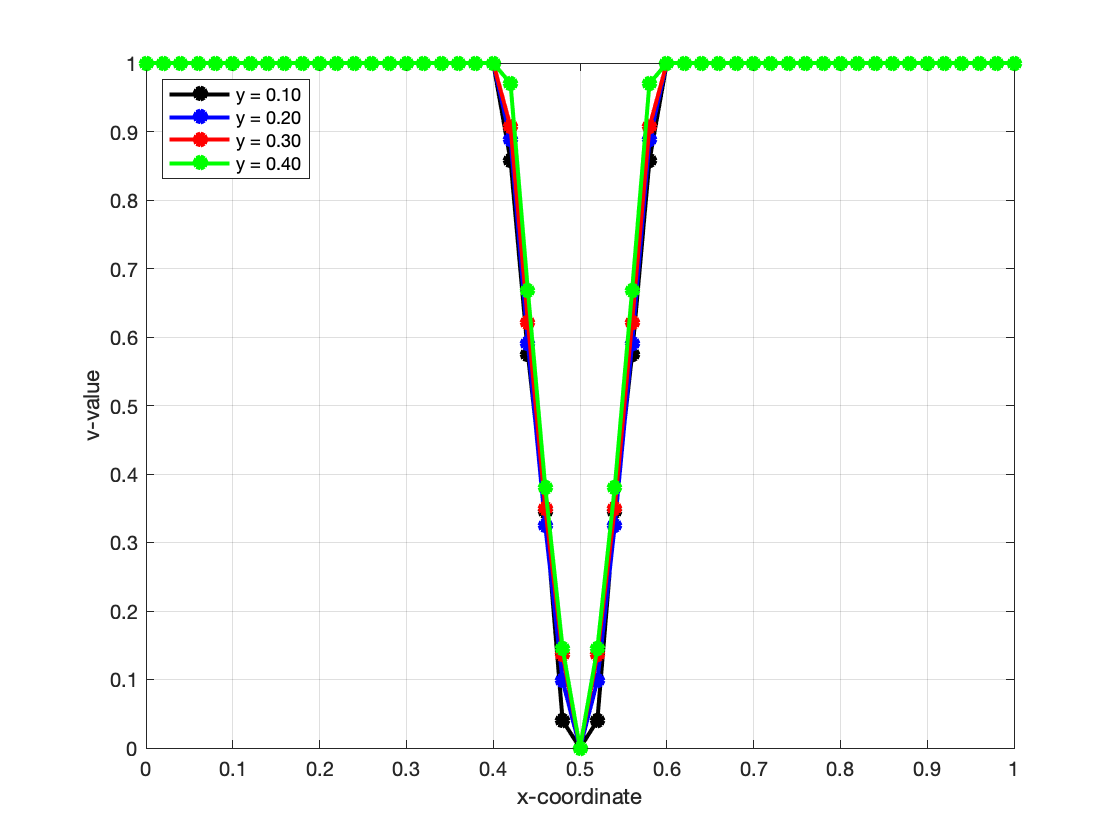}
   \caption{Profiles of the variable $v$ extracted along four distinct lines within the computational domain.}
    \label{fig:v_lines}
\end{figure}

The spatial distribution of the phase field, which represents the fracture pattern, was investigated by analyzing its profile at various distances from the primary crack plane. Figure~\ref{fig:v_lines} displays these profiles plotted along the domain's four distinct horizontal lines: $ 0 \leq x \leq 1$ and  $y=0.1, 0.2, 0.3,$ and $0.4$.  The results demonstrate that the transition zone of the phase field widens with increasing distance from the main fracture at $y=0.1$ (close to the bottom boundary). The sharpest gradient in $v$ is observed at $y=0.4$ (close to the initial crack-tip region), indicating a highly localized crack. In contrast, a slightly broader profile at $y=0.1$ reflects the dissipation of the damage field away from the crack's centerline.

\section{Conclusion}\label{conclusions}

We have introduced a novel phase-field framework that advances fracture simulation through a variational approach to mesh adaptivity. The core of our methodology is a generalization of the classical Ambrosio-Tortorelli (\textsf{AT1}-type) functional into a three-field formulation. This formulation dynamically determines a spatially varying regularization parameter, $\xi(x)$, concurrently with the displacement and phase fields. As $\xi(x)$ dictates the local fracture length scale, its role in the global energy minimization provides a rigorous basis for an adaptive meshing strategy that couples element size to physical behavior.

This model was implemented in a custom solver developed by the authors using the open-source \textsf{deal.II} library and applied to the challenging problem of quasi-static crack propagation under anti-plane shear. Our numerical results demonstrate that the formulation effectively drives local mesh refinement by reducing the length scale parameter in critical regions, resulting in a significant increase in both computational efficiency and approximation accuracy. This enables high-fidelity solutions with substantially fewer degrees of freedom than traditional approaches, which require globally fine meshes. The ability to concentrate computational resources makes our method a tractable and powerful tool, paving the way for future work in large-scale 3D simulations, dynamic fracture, and other multi-physics applications where efficiency is paramount.

\section*{Acknoledgement}
This material is based on work supported by the National Science Foundation under Grant No. 2316905.

\bibliographystyle{plain}
\bibliography{references}

\end{document}